\documentclass{article}
\usepackage[english]{babel}
\usepackage{amsmath, amsfonts, amsthm}
\usepackage{amssymb,mathrsfs, verbatim}
\usepackage{color}

\newtheorem{theorem}{Theorem}[section]

\newtheorem{lemma}[theorem]{Lemma}
\newtheorem{corollary}[theorem]{Corollary}
\newtheorem{definition}[theorem]{Definition}

\begin{document}

\title{On a generalized  Muskat type problem}
\author{Nikolai Chemetov$^1$, Wladimir Neves$^2$}
\date{}
\maketitle

\begin{abstract}

We show the solvability of a multidimensional  Muskat type initial boundary
value problem. 
The proposed mathematical model describing the
transport phenomena of non--homogeneous flow in porous media,
relies on a generalized formulation of the Brinkman equation.
\end{abstract}

\footnotetext[1]{%
CMAF/Universidade de Lisboa, Av. Prof. Gama Pinto, 2, 1649--003 Lisboa,
Portugal E--mail: \textsl{chemetov@ptmat.fc.ul.pt}}

\footnotetext[2]{%
Instituto de Matem\'{a}tica, Universidade Federal do Rio de Janeiro, C.P.
68530, Cidade Universit\'{a}ria 21945--970, Rio de Janeiro, Brazil. E--mail: 
\textsl{wladimir@im.ufrj.br}
}

%
%

%


\section{Introduction}

\label{IN} 
The original Muskat problem was proposed in 1934 by 
Muskat \cite{MUSKAT} to study from Darcy's law
the encroachment of water into an oil sand. Due to applications to oil
reservoir, this problem obtains a great practical interest and also, in view
of the mathematical difficulties, the \textit{\ well--posedness} for
the Muskat problem takes attention of many mathematicians. 

In fact, many important results
concerning the Muskat problem were obtained during the last 20 years. Most
of existence and uniqueness results are related with the situation when
there exists only one moving horizontal interface, that separate two
different fluids. Two regimes were found for the Muskat problem: a stable
regime, when this horizontal interface is stable under small deviations and
an unstable one, which is to say, fingering occurs. 
The stable regime could be realized, if
initially a horizontal interface separates the two fluids with a denser
fluid from below and in the presence of gravity force. In the case of the
stable regime we can mention the following results: Yi \cite{yi}, Siegel,
Caflisch, Howison\ \cite{siegel} shown global--in--time existence for initial
data, that is a small perturbation of a flat interface, which separate two
fluids. Ambrose \cite{ambrose}, C\'{o}rdoba A., C\'{o}rdoba D., Gancedo \cite%
{ACDCFG}, \ Escher, Matioc \cite{JEBVM} proved local--in--time existence and
uniqueness of solutions for initial data in the Sobolev spaces. Further
global well--posedness results were established by Constantin, C\'{o}rdoba,
Gancedo, Strain \cite{PCDCFGRMS} for initial data smaller than an explicitly
computable constant. Nevertheless the well--posedness of the Muskat problem
for general initial data is not known. 
It is also interesting to mention the existence and
non--uniqueness results of weak solutions for the Muskat problem
obtained by C\'{o}rdoba, Faraco, Gancedo \ \cite{cor} and Sz\'{e}kelyhidi
Jr.\ \cite{sze}. The number of constructed weak solutions, corresponding to
the same initial data, is \textit{infinite.} Indeed, the weak solutions considered
have low regularity, and in particular, they do not satisfy a standard
criterium of selection of the unique solution: the renormalized criterium
for the density equation (see Lemma \ref{prop22}). The literature concerning
the Muskat problem is really huge nowadays, we do not pretend here to cover
all of them. We address the reader to the references there in the above
cited works, which seems to us very good to have an up to date scenario of
the Muskat problem.

Whilst Darcy's equation continues to occupy a central place in the study of
flow through porous media, it is only valid for a special class of flows. As
in any fluid flow problem, the range of validity of Darcy's law may be
expressed in terms of the Reynolds number $Re$, which is generally defined
in terms of a characteristic length. In particular, for consolidated porous
media $Re$ is expressed in terms of mean porous size, on the other hand for
unconsolidated porous media, it is in terms of grain size. In any case, it
is claimed that Darcy's law is applicable for $Re$ less than $10,$ where the
viscous forces are predominant. Darcy's law may breakdown for many reasons,
for instance when the $Re$ number is bigger than $100$, or applying for
gases at low pressure, or if the mean porous diameter of the medium is
comparable with the mean free path of the gas, etc. The reader is addressed
to some empirical or almost--empirical modifications of Darcy's law studied
by Scheidegger \cite{Scheidegger}. We also refer to the book of Nield, Bejan 
\cite{DANAB}, where an excellent review of different approaches for
modeling flows in the porous media is given.

It is important to observe that, in the original formulation of the Muskat
problem given by the standard Darcy's law equation, the fluids are assumed
to behave ideally with itself, that is to say, the viscosity of the fluids
only takes place in the constitutive relation of the interaction forces, and
do not in the Cauchy stress tensor, see for instance Appendix B in
Rajagopal, Tao \cite{KRRLT}. In this paper, we follow our original idea
established in \cite{NCWN1}, perturbing the Darcy law equation with a
positive viscosity term, that reduce to a commonly known Brinkman law
equation. Then we formulate an initial boundary value problem assuming
Dirichlet boundary data, and in this way, it is shown the solvability of a
generalized Muskat type problem, that is presented with details in the next
sections.

The Brinkman equation has the form 
\begin{equation*}
\frac{\mu }{\kappa }\,\mathbf{v}-\eta \,\Delta \mathbf{v}=-\nabla p,
\end{equation*}%
where $p$ is the pressure, $\mathbf{v}$ is the velocity field, $\kappa $ is
the permeability and $\eta $ is an effective viscosity, which is not
necessary equal to the dynamic viscosity $\mu $, see for instance 
\cite{DANAB}. In fact, in many applications of the above equation,
also in the original formulation due to Brinkman \cite{brin}, it is used $%
\mu $ instead of $\eta $. Advances in homogenization theory made it possible
rigorously to derive Darcy's and Brinkman's laws from Stokes' equations. As
concluded in \cite{al}, there are three different limits depending on the
size of the periodically arranged obstacles, which respectively lead to
Darcy's, Brinkman's, and Stokes' equations as macroscopic (homogenized)
relation. In the same direction, Sahimi \cite{MS} observes that, if the
length scale where the fluctuations in the velocity field $\mathbf{v}$ are
important is much larger than $O(\sqrt{\kappa })$, then the diffusive term
the Brinkman's equation could be negligible and the Darcy's law is
recovered. Moreover, when the porosity of the medium approaches the
identity, it follows that $\eta /\mu \rightarrow 1$, also $\kappa
\rightarrow \infty $, and the Brinkman's equation approaches the Stokes'
equation. In fact, \ the ratio $\eta /\mu $ depends on the geometry, which
is characterized by \ the porosity, ratio of the void space to the total
volume of the medium, and tortuosity, which represents the hindrance to flow
diffusions imposed by local boundaries or local viscosity of the porous
medium.
Therefore, the Brinkman equation could be interpreted as an interpolation
between the macroscopic Darcy's law and the microscopic Stokes' equation. In
the works Auriault \cite{aur}, Nield,
Bejan \cite{DANAB}, Rajagopal \cite{KRRL} a more complete
comparison study of Darcy's and Brinkman's laws have been developed.

One observes that, the application of the Brinkman law is the natural one in  a high porous medium. 
Moreover, it has been recently employed to problems in
automotive  industry, design of thermal insulations, energy saving \cite{ch},
biomedical hydrodynamics studies \cite{kl},  transport flows in carbonate karst reservoirs \cite{KLLS},
  hydraulic fracturing (mining  of shale gas) \cite{mar}, 
ground water contamination, processes of vaporization--condensation \cite{mo}, among others applications.

Finally, the numerical study of Stokes--Brinkman--Darcy type systems (stationary case) has been developed by 
 Girault, Kanschat, Riviere \cite{gir2}, Ingram \cite{i},  
Layton,  Schieweck, Yotov \cite{L1}.  
The vanishing viscous limit for a Brinkman--Darcy type system has been studied by Kelliher, Temam, Wang \cite{k}.
The global well--posedness of
 the Cauchy problem associated to compressible Brinkman flow is shown in  \cite{m}. 
The Cahn--Hilliard--Brinkman system, modeling a diffuse interface of two--phase separation in porous medium, is investigated
 in the unpublished article \cite{gr}.  Some theoretical results related with the problem studied here,
 can be found in the works \cite{A1} and \cite{den},
 where it was studied  the solvability for 
 the Stokes and Navier--Stokes equations, describing the motion of two fluids. 

\section{Multidimensional Muskat type problem}

\label{MP} 

In this section we establish the multidimensional Muskat type
problem. The model presented here, follows the theory of Continuum Mechanics
for Mixtures. The reader interested to more physical explanations is
addressed to Atkin, Craine \cite{RJAREC}, also Rajagopal \cite{KRRL}, Rajagopal, Tao \cite{KRRLT},
and the treatise of Truesdell, Topin \cite{CTRT}. In particular, for the
study of the interaction terms, we address the reader to Massoudi \cite{MM}.

Generally speaking, the Muskat problem is a piston--like displacement of two
immiscible fluids in a (porous) media domain. A porous media shall be
interpreted like a solid with holes (small ones, mostly interconnected) in
it. Then, the fundamental statement of the Muskat problem is the existence
of two different fluids in a domain $\Omega \subset \mathbb{R}^{d}$ (here $%
\Omega $ is assumed bounded), which immiscible fluids are separated by an
unknown (free) surface $S$ (connected, non--necessarily regular) of
co--dimension one. Moreover, this surface $S$ is not constant in time, and the
main issue is to know if there exists an evolution in time of it.

Since the main motivation is related to planning operation of oil
wells, which is to say, 
the interaction between water and oil, it will be used the subscripts $o,w$,
to distinguish between each phase of the mixture. 
As usual we assume that the mixtures of each fluid and the porous medium is
sufficiently dense, thus for each time--space point $(t,\mathbf{x})$, the
properties of the mixture is well established, for instance the concept of
seepage velocity of each fluid in the porous media.

\medskip Under the assumption that the fluids are immiscible, that is to
say, remain separated for all times during the process of motion
(piston--like), for each $t\geqslant 0$, the bounded smooth domain $\Omega $
is given by the union of two disjoint sets $\Omega _{o}(t)$ and $\Omega
_{w}(t)$, with a common surface $S(t)$ of co--dimension one, where the
surface is assumed Lipschitz. Then, we consider for each $t\geqslant 0$ and $%
\mathbf{x}\in \Omega $ the mixture density $\rho (t,\mathbf{x})$, given by 
\begin{align}
\rho (t,\mathbf{x})= \rho_w(t,\mathbf{x}) \, 1_{\Omega_w(t)}(\mathbf{x}) +
\rho_o(t,\mathbf{x}) \, 1_{\Omega_o(t)}(\mathbf{x}).  \label{RHO}
\end{align}
Moreover, the velocity $\mathbf{v}(t,\mathbf{x})$ of the mixture is 
\begin{equation*}
\mathbf{v}(t,\mathbf{x})=\mathbf{v}_{w}(t,\mathbf{x})\;1_{\Omega _{w}(t)}(%
\mathbf{x})+\mathbf{v}_{o}(t,\mathbf{x})\;1_{\Omega _{o}(t)}(\mathbf{x}),
\end{equation*}%
where $\mathbf{v}_{\iota }(t,\mathbf{x})$, $(\iota =o,w)$, is the velocity
field of each component (also called phase) and is obtained from an average
of the flow rate of the $\iota ^{th}$--phase divided by an unitary area.

\medskip The balance of mass (usually called continuity equation in fluid
dynamics) is given in distribution sense by 
\begin{equation}
\partial _{t}\rho (t,\mathbf{x})+\mathrm{div} \big(\rho (t,\mathbf{x})\;%
\mathbf{v}(t,\mathbf{x})\big)=0.  \label{CEQDF}
\end{equation}

Let us interpret at this point the evolution of the interface $S(t)$. Let $%
\mathbf{w}$ be the velocity of the interface $S(t)$. From the continuity
equation \eqref{CEQDF}, we have on $S(t)$ 
\begin{equation*}
\rho _{o}(\mathbf{v}_{o}-\mathbf{w})\cdot \mathbf{n}_{S}=\rho _{w}(\mathbf{v}%
_{w}-\mathbf{w})\cdot \mathbf{n}_{S}=:m,
\end{equation*}%
where $\mathbf{n}_{S}$ is the unitary normal field to $S(t).$ Since we
assume that the fluids are immiscible, we have $m=0$, i.e. the fluid's flow
does not pass through the interface $S$, hence 
\begin{equation*}
\mathbf{w}\cdot \mathbf{n}_{S}=\mathbf{v}_{o}\cdot \mathbf{n}_{S}=\mathbf{v}%
_{w}\cdot \mathbf{n}_{S}.
\end{equation*}%
This equality describes the evolution of the interface $S(t)$. From the
continuity of the normal velocities of the fluids on $S(t)$ ( $\mathbf{v}%
_{o}\cdot \mathbf{n}_{S}=\mathbf{v}_{w}\cdot \mathbf{n}_{S}$ ) and the
assumption that the fluids are incompressible, that is 
\begin{equation*}
\mathrm{div} \, \mathbf{v}_{\iota }=0\quad \text{in}\ \Omega _{\iota
}(t)\quad \text{for \ }\iota \in M\equiv \{o,w\},
\end{equation*}%
it follows, in the sense of distributions, that 
\begin{equation}
\mathrm{div} \, \mathbf{v}(t,\mathbf{x})=0.  \label{di}
\end{equation}

Now, let us consider the Conservation of Linear Momentum. \ One remarks
that, we are not considering that the fluids behave ideally with itself,
like standard Darcy's law assumption. Consequently, we shall consider an
effective viscosity term in the constitutive relation to the Cauchy stress
tensor. Therefore, we consider the following constitutive relation of the
Cauchy stress tensor for incompressible fluids 
\begin{equation}
\mathbf{T}_{\iota }=-p_{\iota }\;\mathbf{I_{d}}+2\;\eta _{\iota }\;\mathbf{Dv%
}_{\iota },  \label{eq}
\end{equation}%
where $p_{\iota }\geqslant 0,\eta _{\iota }>0$ are respectively the scalar
function called pressure and the effective viscosity of the $\iota ^{th}$%
--phase. Usually the effective viscosity $\eta= \mu / \phi(\mathbf{x})$,
where $\phi > 0$ is the porosity of the medium and $\mu > 0$ is the dynamic
viscosity. Here for simplicity of exposition, we take $\eta= \mu$. Moreover,
one observes that, $\mathbf{I_{d}}$ is the identity map on $\mathbb{R}^{d}$,
and $\mathbf{Dv}_{\iota }$ is the stretching tensor, given in components by 
\begin{equation*}
\lbrack \mathbf{Dv}_{\iota }]_{ij}=\frac{1}{2}\Big(\partial _{\mathbf{x}%
_{j}}[\mathbf{v}_{\iota }]_{i}+\partial _{\mathbf{x}_{i}}[\mathbf{v}_{\iota
}]_{j}\Big).
\end{equation*}%
that is, $\mathbf{D}$ is just the symmetric part of the gradient operator,
instead of the standard notation in Continuum Physics, where $\mathbf{D}$ is
the symmetric part of the gradient operator applied to velocity field.
The constitutive relation for the interaction force terms of the
fluids with the porous media, is given by a functional dependence in the
difference of the velocity of the fluid and the medium. As usual, we assume
linear dependence, and since the porous media velocity is zero, the
interaction force is given by 
\begin{equation}
\mathbf{f}_{\iota }=-h_{\iota }\;\mathbf{v}_{\iota },  \label{fi}
\end{equation}%
with a given nonnegative scalar function $h_{\iota }=h(t,\mathbf{x},\mu
_{\iota })$ (which takes in account the properties of the porous medium).
For instance, we have 
\begin{equation*}
h_{\iota }\big(t,\mathbf{x},\mu _{\iota }(t,\mathbf{x})\big)=\frac{\mu
_{\iota }(t,\mathbf{x})}{k(\mathbf{x})},\qquad \text{($k(\mathbf{x})$ is the
permeability)}.
\end{equation*}%
Certainly, more general conditions on $h_{\iota }$ may be considered. For
instance, we may assume that $h_{\iota }$ is a (positive defined) second
order tensor, which takes in account the anisotropy of the medium. For
simplicity, we just consider here the scalar case. We observe also that, $%
\mu _{\iota }(t,\mathbf{x})=\rho _{\iota }(t,\mathbf{x})\,\nu _{\iota }(t,%
\mathbf{x})$, where $\nu _{\iota }$ is the kinematic viscosity of the $\iota
^{th}$--phase. Once we denote by $\mathbf{g}$ the body force density, 
the generalized Darcy's law is given by 
\begin{equation}
\mathbf{f}_{\iota }(t,\mathbf{x})+\mathrm{div}\,\mathbf{T}_{\iota }(t,%
\mathbf{x})+\rho _{\iota }(t,\mathbf{x})\,\mathbf{g}(t,\mathbf{x})=0,
\label{eq1}
\end{equation}%
that has been proposed by Brinkman \cite{brin} as a new phenomenological
relation between the velocity and the pressure gradient for flows in highly
porous media.

\medskip Due to the assumption of viscous motion of the fluid mixture it is
natural to consider that the velocity $\mathbf{v}$ is continuous on the
interface, that is 
\begin{equation}
\mathbf{v}_{o}=\mathbf{v}_{w} \quad \text{on $S(t)$}.  \label{eq11}
\end{equation}%
Finally, we admit the dynamic boundary condition on $S(t),$ that requires
the continuity of the normal component of the Cauchy stress tensor on $S(t)$:%
\begin{equation}
\mathbf{T}_{o}\mathbf{n}_{S}=\mathbf{T}_{w}\mathbf{n}_{S}.  \label{eq111}
\end{equation}
One observes that, equation \eqref{eq111} gives on the interface the
difference of the oil and water pressures in terms of the viscosity
components. Therefore, could be considered as a corrected physical
interpretation of the Laplace's formula (experimental one), which gives the
difference of the pressures on the interface in terms of the capillarity
pressure.

\medskip From the above considerations \eqref{eq}--\eqref{eq111} we conclude
that the evolution of the mixture velocity $\mathbf{v}$ is described in
distribution sense by 
\begin{align}
& h(t,\mathbf{x},\mu )\;\mathbf{v}(t,\mathbf{x})-\mathrm{div}\,\big(\mu \,%
\mathbf{Dv}(t,\mathbf{x})\big)=-{\nabla }p(t,\mathbf{x})+\rho (t,\mathbf{x}%
)\,\mathbf{g}(t,\mathbf{x}),  \notag \\
& \mu (t,\mathbf{x})=\rho (t,\mathbf{x})\,\nu (t,\mathbf{x}),  \label{TMPV}
\end{align}%
where the kinematic viscosity of the fluid mixture 
\begin{equation*}
\nu (t,\mathbf{x})=\nu _{w}(t,\mathbf{x})\,1_{\Omega _{w}(t)}(\mathbf{x}%
)+\nu _{o}(t,\mathbf{x})\,1_{\Omega _{o}(t)}(\mathbf{x})
\end{equation*}%
is governed, analogously to the density, by the transport equation 
\begin{equation}
\partial _{t}\nu (t,\mathbf{x})+\mathrm{div}\big(\nu (t,\mathbf{x})\;\mathbf{%
v}(t,\mathbf{x})\big)=0  \label{VEQDF}
\end{equation}%
in the distribution sense. The functions $\mu \equiv \rho \nu $ and $p$ are
called the dynamic viscosity and the pressure of the fluid mixture, which
are respectively positive and nonnegative functions, defined by a similar
formula to \eqref{RHO}.

\medskip We recall that, the process is assumed isothermal, thus we do not
consider the Conservation of Energy.

\medskip Let $T>0$ be any fixed real number and $\Omega \subset \mathbb{R}%
^{d}$ (with $d=2$ or $3)$ is an open and bounded domain having a $C^{2}-$
smooth boundary $\Gamma .$ We define by $\Omega _{T}:=(0,T)\times \Omega
,\quad \Gamma _{T}:=(0,T)\times \Gamma .$ Moreover, the outside unitary
normal to $\Omega $ at $\mathbf{x}\in \Gamma $ is denoted by $\mathbf{n}=%
\mathbf{n}(\mathbf{x}).$\ Then, from equations \eqref{CEQDF}, \eqref{di}, %
\eqref{TMPV} and \eqref{VEQDF}, we formulate a generalized Muskat type
initial--boundary value problem, denoted $\mathbf{GMP}$:

For all $(t,\mathbf{x})\in \Omega _{T}$, find $(\rho (t,\mathbf{x}),\nu (t,%
\mathbf{x}),\mathbf{v}(t,\mathbf{x}))$ solution of 
\begin{equation*}
\left\{ \begin{aligned} &\partial_t \rho + \mathbf{v} \cdot{\nabla} \, \rho=
0, \qquad\partial_t \nu + \mathbf{v} \cdot{\nabla}\, \nu= 0, \\[5pt]
&h(t,\mathbf{x}, \rho \nu) \; \mathbf{v} - {\rm div} \big(\rho \nu \,
\,\mathbf{D} \mathbf{v} \big)= -{\nabla} p + \rho \, \mathbf{G}, \quad {\rm
div}\, \mathbf{v}= 0, \\[5pt] &\rho|_{t=0}= \rho_0, \quad
\rho|_{\Gamma_T^-}= \rho_b, \quad \nu|_{t=0}= \nu_0, \quad
\nu|_{\Gamma_T^-}= \nu_b, \quad \mathbf{v}|_{\Gamma_T}= \mathbf{b},
\end{aligned}\right.
\end{equation*}%
where $\mathbf{g}$ is a given vector function, also $\rho _{0}$, $\rho _{b}$%
, $\nu _{0}$, $\nu _{b}$ are given initial--boundary data for the density and
effective viscosity respectively, $\mathbf{b}$ is the boundary data for the
velocity field $\mathbf{v}$ and 
\begin{eqnarray}
\Gamma _{T}^{-} &:&=\left\{ (t,\mathbf{r})\in \Gamma _{T}:(\mathbf{b}\cdot 
\mathbf{n})(t,\mathbf{r})<0\right\} ,  \notag \\
\Gamma _{T}^{+} &:&=\left\{ (t,\mathbf{r})\in \Gamma _{T}:(\mathbf{b}\cdot 
\mathbf{n})(t,\mathbf{r})>0\right\}  \label{flux}
\end{eqnarray}%
called respectively the in--flux, and out--flux boundary zones of the
"oil--water" mixture.

In the following sections we show that, the $\mathbf{GMP}$ is solvable and
describes the motion of immiscible fluids. One of the main difficulties, to
prove the solvability for $\mathbf{GMP}$, is to show the strong convergence
of an approximating sequence for the density function. Another important
issue is the trace of the density function, once it is just assumed
measurable and bounded, in fact, we follow the technical and important
results proved by Boyer in \cite{boyer}. Similar remarks are also posed for
the kinematic viscosity.


\section{Functional notation and auxiliary results}

\label{appendix} 

In this section we present the notations, the definitions of functional
spaces and some useful results, used through the paper.

We will use the standard notations for the Lebesgue function space $%
L^{p}(\Omega ),$ the Sobolev spaces $W^{s,p}(\Omega )$ and $H^{s}(\Omega
)\equiv W^{s,2}(\Omega )$, where a real $s\geqslant 0$ is the smoothness
indices and a real $p\geqslant 1$ is the integrability indices. The vector
counterparts of these spaces are denoted by $\mathbf{L}^{2}(\Omega
)=(L^{2}(\Omega ))^{d}$\ and $\mathbf{H}^{s}(\Omega ):=(H^{s}(\Omega ))^{d}.$

Let us consider the space of functions of $L^{2}-$bounded deformation 
\begin{equation*}
\mathbf{LD}_{0}^{2}(\Omega )=\left\{ \boldsymbol{\psi }\in \mathbf{L}%
^{2}(\Omega ):\;\mathbf{D}\boldsymbol{\psi }\in (\mathbf{L}^{2}(\Omega
))^{d},\;\;\boldsymbol{\psi }=0\;\;\mbox{on}\ \ \Gamma \right\} ,
\end{equation*}%
endowed by the norm$\;||\boldsymbol{\psi }||_{\mathbf{LD}_{0}^{2}(\Omega
)}=||\boldsymbol{\psi }||_{\mathbf{L}^{2}(\Omega )}+||\mathbf{D}\boldsymbol{%
\psi }||_{\mathbf{L}^{2}(\Omega )}.$ Due to the Poincar\'{e} and Korn
inequalities, that is 
\begin{equation}
||\boldsymbol{\psi }||_{\mathbf{H}^{1}(\Omega )}\leqslant C||\mathbf{D}%
\boldsymbol{\psi }||_{\mathbf{L}^{2}(\Omega )},\qquad \forall \boldsymbol{%
\psi }\in \mathbf{LD}_{0}^{2}(\Omega ),  \label{ldd}
\end{equation}%
the space $\mathbf{LD}_{0}^{2}(\Omega )$ shall be considered equivalent to $%
\mathbf{H}_{0}^{1}(\Omega ).$ We also consider the following spaces 
\begin{eqnarray*}
&&\mathbf{V}^{1}(\Omega ):=\Big\{\boldsymbol{\psi }\in \mathbf{H}^{1}(\Omega
):\quad \mathrm{div}\boldsymbol{\psi }=0\quad \text{in }\mathcal{D}%
^{^{\prime }}(\Omega )\Big\}, \\
&&\;\mathbf{V}^{1/2}(\Gamma ):=\Big\{\boldsymbol{\psi }\in \mathbf{H}%
^{1/2}(\Gamma ):\quad \int_{\Gamma }\boldsymbol{\psi }\cdot \mathbf{n}\ d%
\mathbf{x}=0\Big\}.
\end{eqnarray*}

\medskip Let us formulate auxiliary results (Lemmas \ref{prop2}--\ref{prop22}
and Corollary \ref{prop23}), used to prove the main result of our article
(Theorem \ref{theorem}). First, we consider the Stokes type system 
\begin{equation}
\begin{cases}
-\mathrm{div}(\nu \mathbf{Dv})+k(\mathbf{x})\mathbf{v}=-\nabla p+\mathbf{g}%
,\qquad \mathrm{div}\mathbf{v}=0\quad \text{in }\Omega , \\ 
\mathbf{v}=\mathbf{b}\quad \text{on }\Gamma .%
\end{cases}
\label{eq2}
\end{equation}%
The solvability of \eqref{eq2} in the particular case, when $\nu $ is a
constant, $k=0$ \ and $\mathbf{b}=\mathbf{0}$, has been done by Cattabriga 
\cite{cat}. Let us show the existence result for system \eqref{eq2} in a
general situation.

\begin{lemma}
\label{prop2} Let us assume that%
\begin{eqnarray}
\mathbf{b} &\in &\mathbf{V}^{1/2}(\Gamma )\quad \text{and }  \notag \\
\mathbf{g} &\in &\mathbf{L}^{q}(\Omega )\quad \text{for some }q>1\quad \text{%
if }\ d=2;\quad q=\frac{6}{5}\quad \text{if }d=3.  \label{nu0}
\end{eqnarray}%
Let $k(\mathbf{x}),\nu (\mathbf{x})$ be measurable positive functions,
satisfying 
\begin{equation*}
k\in L^{s}(\Omega )\quad \text{for some }s>1\quad \text{if \ }d=2;\quad
s=3/2\quad \text{if }d=3.
\end{equation*}%
\begin{equation}
\nu (\mathbf{x})\in \lbrack \overline{\nu }^{0},\overline{\nu }^{1}]\quad 
\text{a.e. in }\Omega  \label{nu}
\end{equation}%
for some positive real numbers $\overline{\nu }^{0}<\overline{\nu }^{1}.$
Then, there exists a unique weak solution $\mathbf{v}\in \mathbf{V}%
^{1}(\Omega )$ of \eqref{eq2}, satisfying 
\begin{equation}
\Vert \mathbf{v}\Vert _{\mathbf{V}^{1}(\Omega )}\leqslant C(||\mathbf{b}||_{%
\mathbf{V}^{1/2}(\Gamma )}+\Vert \mathbf{g}\Vert _{\mathbf{L}^{q}(\Omega )})
\label{n1}
\end{equation}%
with a positive constant $C$ depending only on $\Omega $ and the data $k,\nu
.$
\end{lemma}

\textbf{Proof.} By Lemma 2.2, p. 24 in \cite{gir}, there exists an operator 
\begin{equation*}
\mathbf{a}\in \mathbf{V}^{1/2}(\Gamma )\longmapsto \mathbf{v}_{a}\in \mathbf{%
V}^{1}(\Omega ),
\end{equation*}%
such that $\mathbf{a}$ is the trace of $\mathbf{v}_{a}$ on $\Gamma $ and
this operator is \textit{linear} and \textit{bounded}, that is 
\begin{equation}
\Vert \mathbf{v}_{a}\Vert _{\mathbf{V}^{1}(\Omega )}\leqslant C\Vert \mathbf{%
a}\Vert _{\mathbf{V}^{1/2}(\Gamma )}  \label{nu1}
\end{equation}%
with $C$ depending only on $\Omega .$ See, also Lemma 3.2, p. 41, in \cite%
{gir}, where the concrete operator is proposed, satisfying \eqref{nu1}.

Let $\mathbf{z}:=\mathbf{v}-\mathbf{v}_{b}$ be the solution of the system%
\begin{equation}
\begin{cases}
-\mathrm{div}(\nu \mathbf{Dz})+k(\mathbf{x})\,\mathbf{z}=-\nabla p+\mathbf{f}%
,\quad \mathrm{div}\mathbf{z}=0\quad \text{in }\Omega , \\ 
\mathbf{z}=\mathbf{0}\quad \text{on }\Gamma%
\end{cases}
\label{eq22}
\end{equation}%
in the sense of the integral identity\textbf{\ \ }%
\begin{equation*}
\int_{\Omega }\left\{ \nu \,\mathbf{Dz}:\mathbf{D}\boldsymbol{\psi }+k\,%
\mathbf{z}\boldsymbol{\psi }\right\} \ d\mathbf{x}=-\int_{\Omega }\left\{
\nu \,\mathbf{Dv}_{b}:\mathbf{D}\boldsymbol{\psi }\;+k\mathbf{v}_{b}%
\boldsymbol{\psi }-\mathbf{g}\boldsymbol{\psi }\right\} \ d\mathbf{x}
\end{equation*}%
for any $\ \boldsymbol{\psi }\in \mathbf{V}^{1}(\Omega ),$ such that $%
\boldsymbol{\psi }=0$\ on $\Gamma .$ Here $\mathbf{f}:=\mathrm{div}(\nu 
\mathbf{Dv}_{b})-k\mathbf{v}_{b}+\mathbf{g}.$\ \ The solvability of %
\eqref{eq22} follows from the Lax--Milgram theorem (see Proposition 2.2, p.12
in \cite{const}), the Korn inequality \eqref{ldd}, assumptions \eqref{nu0}--%
\eqref{nu} and the embedding theorem 
\begin{equation}
W^{1,2}(\Omega )\subset L^{r}(\Omega )\ \ \text{for }r<\infty \quad \text{
if }\ d=2;\quad r=6\quad \text{if }d=3.  \label{s}
\end{equation}%
By a standard argument, we show that, the function $\mathbf{v}:=\mathbf{z}+%
\mathbf{v}_{b}$ \ is the unique weak solution of \eqref{eq2}.

Due to \eqref{ldd} the function$\ \mathbf{z}$ satisfies the following
estimates%
\begin{align*}
\overline{\nu }^{0}\Vert \mathbf{Dz}\Vert _{\mathbf{L}^{2}(\Omega )}^{2}+&
\int_{\Omega }k|\mathbf{z}|^{2}\,d\mathbf{x}\leqslant \int_{\Omega }\nu |%
\mathbf{Dz}|^{2}\,d\mathbf{x}+\int_{\Omega }k|\mathbf{z}|^{2}\,d\mathbf{x} \\
\leqslant & \int_{\Omega }\nu |\mathbf{Dz}:\mathbf{D\mathbf{v}_{b}}|\,d%
\mathbf{x}+\int_{\Omega }k|\mathbf{z}||\mathbf{v}_{b}|\,d\mathbf{x}%
+\int_{\Omega }|\mathbf{z}||\mathbf{g}|d\mathbf{x} \\
& \leqslant \frac{\overline{\nu }^{0}}{2}\Vert \mathbf{Dz}\Vert _{\mathbf{L}%
^{2}(\Omega )}^{2}+\int_{\Omega }k|\mathbf{z}|^{2}\,d\mathbf{x+}C||\mathbf{v}%
_{b}||_{\mathbf{V}^{1}(\Omega )}^{2} \\
& +C\int_{\Omega }k|\mathbf{v}_{b}|\,^{2}d\mathbf{x}+C\Vert \mathbf{g}\Vert
_{\mathbf{L}^{q}(\Omega )}^{2},
\end{align*}%
where the constant $C$ depends only on $\Omega $ and $\overline{\nu }^{0},%
\overline{\nu }^{1}.$\ By \eqref{nu0}--\eqref{nu}, \eqref{nu1} \ and \eqref{s}%
, we have 
\begin{equation*}
\Vert \mathbf{Dz}\Vert _{\mathbf{L}^{2}(\Omega )}\leqslant C(||\mathbf{b}||_{%
\mathbf{V}^{1/2}(\Gamma )}+\Vert \mathbf{g}\Vert _{\mathbf{L}^{q}(\Omega )}),
\end{equation*}%
that implies \eqref{n1}.$\hfill \;\blacksquare $

\bigskip

Next we give a time dependent generalization of Lemma \ref{prop2}.

\begin{lemma}
\label{trace} Let us assume that%
\begin{eqnarray*}
\mathbf{b} &\in &L^{2}(0,T;\mathbf{V}^{1/2}(\Gamma ))\quad \text{and} \\
\mathbf{g} &\in &L^{2}(0,T;\mathbf{L}^{q}(\Omega ))\quad \text{for }q>1\quad 
\text{if }\ d=2;\quad q=\frac{6}{5}\quad \text{if }d=3.
\end{eqnarray*}%
Let $k(t,\mathbf{x}),\nu (t,\mathbf{x})$ be measurable positive functions and%
\begin{equation*}
k\in L^{\infty }(0,T;L^{s}(\Omega ))\quad \text{for }s>1\quad \text{if \ }%
d=2;\quad s=3/2\quad \text{if }d=3,
\end{equation*}%
\begin{equation*}
\nu \in \lbrack \overline{\nu }^{0},\overline{\nu }^{1}]\quad \text{a.e. in }%
\Omega _{T}
\end{equation*}%
for some positive real numbers $\overline{\nu }^{0}<\overline{\nu }^{1}.$
Then, there exists a unique weak solution $\mathbf{v}\in L^{2}(0,T;\mathbf{V}%
^{1}(\Omega ))$ of \eqref{eq2} for a.e. $t\in (0,T),$ satisfying%
\begin{equation}
\Vert \mathbf{v}\Vert _{L^{2}(0,T;\mathbf{V}^{1}(\Omega ))}\leqslant C
\label{ca00}
\end{equation}%
with a positive constant $C$ depending only on $\Omega $ and the data $k,\nu
,\mathbf{b},\mathbf{g}.$
\end{lemma}

\textbf{Proof.} Since $\mathbf{b}\in L^{2}(0,T;\mathbf{V}^{1/2}(\Gamma ))$,
then there exists a sequence of \textit{simple} functions $\{\mathbf{b}_{n}=%
\mathbf{b}_{n}(t)\in \mathbf{V}^{1/2}(\Gamma ):t\in (0,T)\}_{n=1}^{\infty },$
such that%
\begin{eqnarray}
||\mathbf{b}_{n}(t)-\mathbf{b}(t)||_{\mathbf{V}^{1/2}(\Gamma )} &\rightarrow
&0\quad \quad \text{for a.e. }t\in \lbrack 0,T],  \notag \\
\Vert \mathbf{b}_{n}-\mathbf{b}\Vert _{L^{2}(0,T;\mathbf{V}^{1/2}(\Gamma ))}
&\rightarrow &0,  \label{ca0}
\end{eqnarray}%
in accordance with the definitions, given in \cite{EG}, p. 649--650 and the
proof, given in Theorem 1, p. 133 \cite{y} (see, \textit{The "if" \ part}).
Applying \eqref{nu1}, there exists a sequence of \textit{simple} functions $%
\{\mathbf{v}_{b,n}(t)\in \mathbf{V}^{1}(\Omega ):t\in (0,T)\}_{n=1}^{\infty
},$ such that $\mathbf{v}_{b,n}(t)=\mathbf{b}_{n}(t)$ on $\Gamma $ for $t\in
(0,T),$ satisfying%
\begin{eqnarray}
\Vert \mathbf{v}_{b,n}(t)\Vert _{\mathbf{V}^{1}(\Omega )} &\leqslant &C\Vert 
\mathbf{b}_{n}(t)\Vert _{\mathbf{V}^{1/2}(\Gamma )},  \notag \\
\Vert \mathbf{v}_{b,n}(t)-\mathbf{v}_{b,n^{\prime }}(t)\Vert _{\mathbf{V}%
^{1}(\Omega )} &\leqslant &C\Vert \mathbf{b}_{n}(t)-\mathbf{b}_{n^{\prime
}}(t)\Vert _{\mathbf{V}^{1/2}(\Gamma )}  \label{ca}
\end{eqnarray}%
for any $\ \forall n,n^{\prime }=1,2,....$ Here $C$ is a positive constant
depending only on $\Omega .$ Therefore the sequence of \textit{simple}
functions $\{\mathbf{v}_{b,n}\}_{n=1}^{\infty }$ is the Cauchy sequence in $%
L^{2}(0,T;\mathbf{V}^{1}(\Omega ))$ by \eqref{ca0}-\eqref{ca}. Since $%
L^{2}(0,T;\mathbf{V}^{1}(\Omega ))$ is a Banach space, there exists a 
\textit{measurable} $\mathbf{v}_{b}\in L^{2}(0,T;\mathbf{V}^{1}(\Omega ))$,
such that $\mathbf{v}_{b}(t)=\mathbf{b}(t)$ on $\Gamma $ for a.e. $t\in
(0,T) $ and%
\begin{equation}
\mathbf{v}_{b,n}\rightarrow \mathbf{v}_{b}\quad \quad \text{in }L^{2}(0,T;%
\mathbf{V}^{1}(\Omega ))\quad \quad \text{as }n\rightarrow \infty .
\label{vb}
\end{equation}%
Let $\mathbf{z}:=\mathbf{v}-\mathbf{v}_{b}$ be the solution of the system%
\begin{equation}
\begin{cases}
-\mathrm{div}(\nu \mathbf{Dz})+k(\mathbf{x})\,\mathbf{z}=-\nabla p+\mathbf{f}%
,\quad \mathrm{div}\mathbf{z}=0\quad \text{in }\Omega _{T}, \\ 
\mathbf{z}(t)=\mathbf{0}\quad \text{on }\Gamma \quad \quad \text{for a.e. }%
t\in \lbrack 0,T],%
\end{cases}
\label{eq202}
\end{equation}%
in the sense of the integral identity%
\begin{equation*}
\int_{\Omega _{T}}\left\{ \nu \,\mathbf{Dz}:\mathbf{D}\boldsymbol{\psi }+k\,%
\mathbf{z}\boldsymbol{\psi }\right\} \ d\mathbf{x}dt=-\int_{\Omega
_{T}}\left\{ \nu \,\mathbf{Dv}_{b}:\mathbf{D}\boldsymbol{\psi }\;+k\,\mathbf{%
v}_{b}\boldsymbol{\psi }-\mathbf{g}\boldsymbol{\psi }\right\} \ d\mathbf{x}dt
\end{equation*}%
for any $\ \boldsymbol{\psi }\in L^{2}(0,T;\mathbf{V}^{1}(\Omega )),$ such
that $\boldsymbol{\psi }(t)=\mathbf{0}$\ on $\Gamma $ for a.e. $t\in (0,T).$%
\ Here $\mathbf{f}:=\mathrm{div}(\nu \,\mathbf{Dv}_{b})-k\,\mathbf{v}_{b}+%
\mathbf{g}.$\ \ Due to the Lax--Milgram theorem there exists a solution $%
\mathbf{z}\in L^{2}(0,T;\mathbf{V}^{1}(\Omega ))$ of \eqref{eq202}.
Obviously, $\mathbf{v}=$ $\mathbf{z}+\mathbf{v}_{b}$ is \ the unique weak
solution of \eqref{eq2} for a.e. $t\in (0,T),$ satisfying \eqref{ca00}, that
can be shown by the same way as \eqref{n1}.$\hfill \;\blacksquare $

\bigskip Let us consider the linear transport equation in bounded domains,
that is 
\begin{equation}
\begin{cases}
\partial _{t}\rho +\mathrm{div}\big(\mathbf{v}\rho \big)=0\quad \text{in }%
\Omega _{T}, \\ 
\rho |_{t=0}=\rho _{0}\quad \text{in }\Omega \text{\qquad and\qquad }\rho
=\rho _{b}\quad \text{on }\Gamma _{T}^{-},%
\end{cases}
\label{transport}
\end{equation}%
where $\Gamma _{T}^{-}$ is defined by \eqref{flux}. The proof of the second
auxiliary result, which we follow here, has been obtained by Boyer \cite%
{boyer}, see Theorem 4.1. 
For convenience, we denote by $d\mathbf{r}dt$ the induced surface measure on 
$\Gamma _{T}$, and also, we introduce the measure $d\mu :=(\mathbf{b}\cdot 
\mathbf{n})\,d\mathbf{x}dt$ on $\Gamma _{T}.$ \ From Jordan's decomposition,
we have $\mu =\mu ^{+}-\mu ^{-}$, where $d\mu ^{\pm }:=(\mathbf{b}\cdot 
\mathbf{n})^{\pm }\,d\mathbf{x}dt.$ \ Then, we have the following

\begin{lemma}
\label{prop22} Assume that $\rho _{0}\in L^{\infty }(\Omega ),$ $\rho
_{b}\in L^{\infty }(\Gamma _{T};\mu ^{-}),$ $\mathbf{v}\in L^{2}(0,T;\mathbf{%
V}^{1}(\Omega ))$ and $\mathbf{b}\in L^{2}(0,T;\mathbf{V}^{1/2}(\Gamma ))$
are given functions, such that $\mathbf{v}|_{\Gamma _{T}}=\mathbf{b.}$ Then,
there exists a unique pair 
\begin{equation*}
(\rho ,\rho ^{o})\in L^{\infty }(\Omega _{T})\times L^{\infty }(\Gamma
_{T};\mu ^{+}),
\end{equation*}%
which is a weak solution of \eqref{transport}, satisfying the integral
equality%
\begin{align}
\iint_{\Omega _{T}}\rho \;\Big(\phi _{t}& +(\mathbf{v}\cdot {\nabla )}\phi %
\Big)\;d\mathbf{x}dt+\int_{\Omega }\rho _{0}\;\phi (0)\;d\mathbf{x}  \notag
\\
& =\iint_{\Gamma _{T}}\rho ^{o}\,\phi \,d\mu ^{+}-\iint_{\Gamma _{T}}\rho
_{b}\,\phi \,d\mu ^{-}  \label{1}
\end{align}%
for each test function $\phi \in C_{c}^{\infty }((-\infty ,T)\times \mathbb{R%
}^{d})$. Moreover, the pair $(\rho ,\rho ^{o})$ satisfies 
\begin{equation*}
\Vert \rho \Vert _{L^{\infty }(\Omega _{T})},\Vert \rho ^{o}\Vert
_{L^{\infty }(\Gamma _{T};\mu ^{+})}\leqslant \max \{\Vert \rho _{0}\Vert
_{L^{\infty }(\Omega )},\Vert \rho _{b}\Vert _{L^{\infty }(\Gamma _{T};\mu
^{-})}\}.
\end{equation*}%
The pair $(\rho ,\rho ^{o})$ is a renormalized solution of \eqref{transport}%
, which means that for any function $\beta \in C^{1}(\mathbb{R})$, the pair 
\begin{equation*}
\alpha =(\beta (\rho )\quad \text{in }\Omega _{T},~\beta (\rho ^{o})\quad 
\text{on }\Gamma _{T}^{+})
\end{equation*}%
is the unique solution of the system%
\begin{equation}
\begin{cases}
\partial _{t}\alpha +\mathrm{div}\big(\mathbf{v}\alpha \big)=0\quad \text{in 
}\Omega _{T}\text{\quad and\quad }\alpha =\beta (\rho ^{o})\quad \text{on }%
\Gamma _{T}^{+},\text{\qquad } \\ 
\alpha |_{t=0}=\beta (\rho _{0})\quad \text{in }\Omega \text{\quad and\quad }%
\alpha =\beta (\rho _{b})\quad \text{on }\Gamma _{T}^{-}%
\end{cases}
\label{2}
\end{equation}%
in the distributional sense, i.e. $\alpha $\ satisfies the equality%
\begin{align}
\iint_{\Omega _{T}}\alpha \;\Big(\phi _{t}& +(\mathbf{v}\cdot {\nabla )}\phi %
\Big)\;d\mathbf{x}dt+\int_{\Omega }\beta (\rho _{0})\;\phi (0)\;d\mathbf{x} 
\notag \\
& =\iint_{\Gamma _{T}}\alpha \,\phi \,d\mu ^{+}-\iint_{\Gamma _{T}}\beta
(\rho _{b})\,\phi \,d\mu ^{-}  \label{3}
\end{align}%
for each test function $\phi \in C_{c}^{\infty }((-\infty ,T)\times \mathbb{R%
}^{d})$.
\end{lemma}

\bigskip As a consequence of Lemma \ref{prop22} we have the next result,
which is similar to that one established in \cite{boyer} (see Lemma 4.2) for
approximated solutions. We present here the proof of this result, since the
same argument will be used in the sequel.

\begin{corollary}
\label{prop23} Let the data $\rho _{0},\rho _{b},$ $\mathbf{b,}$ $\mathbf{v}$
fulfill the assumptions of Lemma \ref{prop22}. Moreover we assume that 
\begin{equation*}
\rho _{0}(\mathbf{x})\in \lbrack \overline{\rho }^{0},\overline{\rho }%
^{1}]\quad \text{a.e. in }\Omega ,\qquad \rho _{b}(t,\mathbf{x})\in \lbrack 
\overline{\rho }^{0},\overline{\rho }^{1}]\quad \text{a.e. in }\Gamma
_{T}^{-}
\end{equation*}%
for some real numbers $\overline{\rho }^{0}<\overline{\rho }^{1}.$ Then the
unique weak solution $(\rho ,\rho ^{o})\in L^{\infty }(\Omega _{T})\times
L^{\infty }(\Gamma _{T};\mu ^{+})$\ of \ \eqref{transport} satisfies 
\begin{equation}
\rho \in \lbrack \overline{\rho }^{0},\overline{\rho }^{1}]\quad \text{a.e.
in }\Omega _{T}\quad \text{and\quad }\rho ^{o}\in \lbrack \overline{\rho }%
^{0},\overline{\rho }^{1}]\quad \text{a.e. in }\Gamma _{T}^{+}.  \label{uff}
\end{equation}
\end{corollary}

\textbf{Proof. \ } In view of Lemma \ \ref{prop22}, \ the unique solution $%
(\rho ,\rho ^{o})$ of \eqref{transport} is a renormalized solution, which
means that for any positive function $\beta \in C^{1}(\mathbb{R})$ the pair $%
\alpha =(\beta (\rho ),\beta (\rho ^{o}))$ is the unique weak solution of
system \eqref{2}. In particular if we choose $\beta $ with $\mathrm{supp}%
(\beta )\subset \overline{\mathbb{R} \setminus \lbrack \overline{\rho }^{0},%
\overline{\rho }^{1}]},$ then $\alpha $ solves the system%
\begin{equation*}
\begin{cases}
\partial _{t}\alpha +\mathrm{div}\big(\mathbf{v}\alpha \big)=0\quad \text{in 
}\Omega _{T}\text{\quad and\quad }\alpha =\beta (\rho ^{o})\quad \text{on }%
\Gamma _{T}^{+}, \\ 
\alpha |_{t=0}=\beta (\rho _{0})\equiv 0\quad \text{in }\Omega \text{\quad
and\quad }\alpha =\beta (\rho _{b})\equiv 0\quad \text{on }\Gamma _{T}^{-}.%
\end{cases}%
\end{equation*}%
Since the zero is the unique solution of this system, we conclude that%
\begin{equation*}
\alpha =(\beta (\rho ),\beta (\rho ^{o}))\equiv 0,
\end{equation*}%
which implies \eqref{uff}. $\hfill \;\blacksquare $


\section{Single--Phase Filtration. Solvability of $\mathbf{GMP}$}

\label{SSOP1} 

Let us first consider one non--homogeneous fluid, which motion is described
by the system $\mathbf{GMP,}$ and assume that the data $\rho _{0}$, $\rho
_{b},\ \nu _{0},\nu _{b},\ \mathbf{b}$, $\mathbf{g}$, the function $h$\
satisfy the following properties 
\begin{eqnarray}
\mathbf{b} &\in &L^{2}(0,T;\mathbf{V}^{1/2}(\Gamma )),\qquad  \notag \\
\mathbf{g} &\in &L^{2}(0,T;\mathbf{L}^{q}(\Omega ))\quad \text{for }q>1\quad 
\text{if \ }d=2;\quad q=\frac{6}{5}\quad \text{if }d=3,  \notag \\
\rho _{0}(\mathbf{x}) &\in &[\overline{\rho }^{0},\overline{\rho }^{1}]\quad 
\text{a.e. in }\Omega \text{,\qquad }\rho _{b}(t,\mathbf{x})\in \lbrack 
\overline{\rho }^{0},\overline{\rho }^{1}]\quad \text{a.e. in }\Gamma
_{T}^{-},  \notag \\
\nu _{0}(\mathbf{x}) &\in &[\overline{\nu }^{0},\overline{\nu }^{1}]\quad 
\text{a.e. in }\Omega \text{,\qquad }\nu _{b}(t,\mathbf{x})\in \lbrack 
\overline{\nu }^{0},\overline{\nu }^{1}]\quad \text{a.e. in }\Gamma _{T}^{-}
\label{reg2}
\end{eqnarray}%
for some positive real numbers $\overline{\rho }^{0}\leqslant \overline{\rho 
}^{1},$ $\overline{\nu }^{0}\leqslant \overline{\nu }^{1}$ and, $h$ is a
Carath\'{e}odory function, more precisely for each $r\geqslant 0$, $h(\cdot ,%
\mathbf{\cdot },r)$ is a measurable function, and for almost all $(t,\mathbf{%
x})\in \Omega _{T}$, the map $0\leqslant r\mapsto h(t,\mathbf{x},r)$ is
continuous. Moreover, we assume for almost all $(t,\mathbf{x})\in \Omega
_{T} $, and each $r\geqslant 0$, 
\begin{eqnarray}
&&0\leqslant h(t,\mathbf{x},r)\leqslant h_{0}(t,\mathbf{x})\;r^{m}\quad 
\text{for some }\;\;m\in \mathbb{R}^{+},  \notag \\
&&h_{0}\in L^{\infty }(0,T;L^{s}(\Omega ))\quad \text{for }s>1\quad \text{if
\ }d=2;\quad s=3/2\quad \text{if }d=3.  \label{reg3}
\end{eqnarray}

The following definition tells us in which sense we consider that a triple
is a weak solution to the problem $\mathbf{GMP}$.

\begin{definition}
\label{DSSMP} The triple $((\rho ,\rho ^{o}),~(\nu ,\nu ^{o}),~\mathbf{v})$
is called a weak solution to the problem $\mathbf{GMP}$, if $(\rho ,\rho
^{o}),$\ $(\nu ,\nu ^{o})\in L^{\infty }(\Omega _{T})\times L^{\infty
}(\Gamma _{T};\mu ^{+}),$ $\ \mathbf{v}\in L^{2}(0,T;\mathbf{V}^{1}(\Omega
)),$ which satisfy the integral identities%
\begin{align*}
\iint_{\Omega _{T}}\rho \;\Big(\phi _{t}& +(\mathbf{v}\cdot {\nabla )}\phi %
\Big)\;d\mathbf{x}dt+\int_{\Omega }\rho _{0}\;\phi (0)\;d\mathbf{x} \\
& =\iint_{\Gamma _{T}}\rho ^{o}\,\phi \,d\mu ^{+}-\iint_{\Gamma _{T}}\rho
_{b}\,\phi \,d\mu ^{-},
\end{align*}%
\begin{align*}
\iint_{\Omega _{T}}\nu \;\Big(\phi _{t}& +(\mathbf{v}\cdot {\nabla )}\phi %
\Big)\;d\mathbf{x}dt+\int_{\Omega }\nu _{0}\;\phi (0)\;d\mathbf{x} \\
& =\iint_{\Gamma _{T}}\nu ^{o}\,\phi \,d\mu ^{+}-\iint_{\Gamma _{T}}\nu
_{b}\,\phi \,d\mu ^{-},
\end{align*}%
\begin{equation*}
\int_{\Omega }\Big(h(t,\mathbf{x},\nu \rho )\;\mathbf{v}\cdot \boldsymbol{%
\psi }+\nu \rho \,\mathbf{Dv}:\mathbf{D}\boldsymbol{\psi }\Big)\;d\mathbf{x}%
=\int_{\Omega }\rho \,\mathbf{g}\cdot \boldsymbol{\psi }\;d\mathbf{x}\quad 
\text{for a.a. }t\in (0,T),
\end{equation*}%
for each test functions $\phi \in C_{c}^{\infty }((-\infty ,T)\times \mathbb{%
R}^{d})$ and $\boldsymbol{\psi }\in \mathbf{V}^{1}(\Omega ),$ such that $%
\boldsymbol{\psi }=0$ on $\Gamma .$ \ Moreover the trace of $\mathbf{v}$ is
equal to $\mathbf{b}$ on $\Gamma _{T}.$
\end{definition}

\bigskip

\begin{theorem}
\label{theorem} Under assumptions \eqref{reg2}, \eqref{reg3} the problem $%
\mathbf{GMP}$ has a weak solution.
\end{theorem}

In the following subsection \ref{apr}, we prove this theorem.


\subsection{Schauder's fixed point argument}

\label{apr} 

\bigskip We apply the Schauder fixed point argument to show the existence
result. To begin, we consider the closed convex subset 
\begin{equation}
\mathcal{Z}=\{(\rho ,\nu )\in L^{2}(\Omega _{T})^{2}:\rho \in \lbrack 
\overline{\rho }^{0},\overline{\rho }^{1}],\quad \nu \in \lbrack \overline{%
\nu }^{0},\overline{\nu }^{1}]\quad \text{a.e. in }\Omega _{T}\}
\label{eqt2.5}
\end{equation}%
of the Banach space $L^{2}(\Omega _{T})^{2}$, with the norm 
\begin{equation*}
||(\rho ,\nu )||_{L^{2}(\Omega _{T})^{2}}\!\!:=||\rho ||_{L^{2}(\Omega
_{T})}\!\!+||\nu ||_{L^{2}(\Omega _{T})}.
\end{equation*}

\bigskip Let $(\overline{\rho },\overline{\nu })$ be an arbitrary fixed
element of $\mathcal{Z}$, and consider the coupled system 
\begin{equation}
\begin{cases}
-\mathrm{div}(\overline{\rho }\ \overline{\nu }\mathbf{Dv})+h(t,\mathbf{x},%
\overline{\rho }\ \overline{\nu })\mathbf{v}=-\nabla p+\overline{\rho }%
\mathbf{g},\quad \mathrm{div}\mathbf{v}=0\quad \text{in }\Omega _{T}, \\ 
\mathbf{v}=\mathbf{b}\quad \text{on }\Gamma _{T},%
\end{cases}
\label{a}
\end{equation}%
\begin{equation}
\begin{cases}
\partial _{t}\rho +\mathrm{div}(\mathbf{v}\rho )=0\quad \text{in }\Omega
_{T}, \\ 
\rho =\rho ^{o}\quad \text{on }\Gamma _{T}^{+}, \\ 
\rho |_{t=0}=\rho _{0}\quad \text{in }\Omega , \\ 
\rho =\rho _{b}\quad \text{on }\Gamma _{T}^{-},%
\end{cases}%
\quad 
\begin{cases}
\partial _{t}\nu +\mathrm{div}(\mathbf{v}\nu )=0\quad \text{in }\Omega _{T},
\\ 
\nu =\nu ^{o}\quad \text{on }\Gamma _{T}^{+}, \\ 
\nu |_{t=0}=\nu _{0}\quad \text{in }\Omega , \\ 
\nu =\nu _{b}\quad \text{on }\Gamma _{T}^{-}.%
\end{cases}
\label{b}
\end{equation}

\bigskip

The solvability result for this system is presented in the following

\begin{lemma}
\label{PPA} For each $(\overline{\rho },\overline{\nu })\in \mathcal{Z}$,
there exists a unique solution $(\rho ,\nu ,\mathbf{v})$ \ of system %
\eqref{a}--\eqref{b}, such that%
\begin{equation}
(\rho ,\nu )\in \mathcal{Z},\qquad \Vert \mathbf{v}\Vert _{L^{2}(0,T;\mathbf{%
V}^{1}(\Omega ))}\leqslant C.  \label{EUU}
\end{equation}%
Hereupon, $C$ is a positive constant depending only on data \eqref{reg2}, %
\eqref{reg3}.
\end{lemma}

\textbf{Proof.} Due to \eqref{reg2}, \eqref{reg3} and Lemma \ref{trace}, the
Stokes type system \eqref{a} has a unique solution $\mathbf{v}=\mathbf{v}(t,%
\mathbf{x})$ in $L^{2}(0,T;\mathbf{V}^{1}(\Omega ))$, such that 
\begin{equation*}
\Vert \mathbf{v}\Vert _{L^{2}(0,T;\mathbf{V}^{1}(\Omega ))}\leqslant C,
\end{equation*}%
where the constant $C$ depends only on $\Omega ,$ $\mathbf{b}$, $\mathbf{g}$%
, $\overline{\nu }^{i},$ $\overline{\rho }^{i},$\ $i=0,1.$\ The last one
permits to apply Lemma \ref{prop22} and Corollary \ref{prop23} with the help
of \eqref{reg2}. Hence systems \eqref{b} have unique solutions $(\rho ,\rho
^{o}),$\ $(\nu ,\nu ^{o})\in L^{\infty }(\Omega _{T})\times L^{\infty
}(\Gamma _{T};\mu ^{+})$, satisfying the estimates%
\begin{eqnarray*}
\rho &\in &[\overline{\rho }^{0},\overline{\rho }^{1}]\quad \text{a.e. in }%
\Omega _{T}\quad \text{and\quad }\rho ^{o}\in \lbrack \overline{\rho }^{0},%
\overline{\rho }^{1}]\quad \text{a.e. in }\Gamma _{T}^{+}, \\
\nu &\in &[\overline{\nu }^{0},\overline{\nu }^{1}]\quad \text{a.e. in }%
\Omega _{T}\quad \text{and\quad }\nu ^{o}\in \lbrack \overline{\nu }^{0},%
\overline{\nu }^{1}]\quad \text{a.e. in }\Gamma _{T}^{+}.
\end{eqnarray*}%
Therefore the solvability of coupled systems \eqref{a}--\eqref{b} is shown
and $(\rho ,\nu ,\mathbf{v})$ \ satisfies \eqref{EUU}. $\hfill
\;\blacksquare $

\bigskip Recall that, solving \eqref{a}--\eqref{b} we have constructed the
operator 
\begin{equation*}
P:\mathcal{Z}\rightarrow \mathcal{Z},\quad \mathbf{(}\rho ,\nu )=P(\overline{%
\rho }\mathrm{,}\overline{\nu }),\;\quad \forall \,(\overline{\rho },%
\overline{\nu })\in \mathcal{Z}.
\end{equation*}%
To find a fixed point of $P$ (by Schauder's theorem), which will be a
solution of the system $\mathbf{GMP}$, it is enough to show that $P(\mathcal{%
Z})$ is a relatively compact subset of the Banach space $L^{2}(\Omega
_{T})^{2}$, and also $P$ is a continuous operator with respect to the norm $%
||(\cdot ,\cdot )||_{L^{2}(\Omega _{T})^{2}}$. \ First, let us show the
following

\begin{lemma}
\label{PPA2} The set $P(\mathcal{Z})$ is relatively compact in $L^{2}(\Omega
_{T})^{2}.$
\end{lemma}

\textbf{Proof.} \ Let us consider an arbitrary sequence $\left\{ (\overline{%
\rho }^{n},\overline{\nu }^{n})\in \mathcal{Z}\right\} _{n=0}^{\infty }$.
For each $n,$ \ due to Lemma \ref{PPA},\ there exists a unique solution $%
(\rho ^{n},\rho ^{o,n}),$ $(\nu ^{n},\nu ^{o,n}),$ $\mathbf{v}^{n}$\ \ of
system \eqref{a}--\eqref{b}, \ which satisfies \eqref{EUU}. We have $(\rho
^{n},\nu ^{n})=P(\overline{\rho }^{n},\overline{\nu }^{n})\in \mathcal{Z}.$
\ Since the set $\mathcal{Z}$ is bounded in $L^{\infty }(\Omega _{T})\subset
L^{2}(\Omega _{T}),$ there exists a suitable subsequence, which we continue
to denote by the same index $"n"$, just for convenience of reading. Then, 
\begin{align}
& \rho ^{n},\;\nu ^{n}\rightharpoonup \rho ,\;\nu \qquad \star \text{-weakly
in }\;L^{\infty }(\Omega _{T})\quad \text{and}\quad \text{weakly in }%
L^{2}(\Omega _{T}),  \notag \\
& \mathbf{v}^{n}\rightharpoonup \mathbf{v}\qquad \text{weakly in }%
\;L^{2}(0,T;\mathbf{V}^{1}(\Omega ))  \label{up}
\end{align}%
with $\rho ,\nu ,\mathbf{v}$, satisfying \eqref{EUU}.

\medskip In the sequel the issue is to prove that, the weak convergences %
\eqref{up} implies the strong convergence in $L^{2}(\Omega _{T})$ for a
suitable subsequence of $\left\{ \rho ^{n}\right\} _{n=0}^{\infty }.$ The
strong convergence for a subsequence of $\left\{ \nu ^{n}\right\}
_{n=0}^{\infty }$\ can be shown analogously. First, there exists a
subsequence of $\left\{ \rho ^{n},(\rho ^{n})^{2}\right\} _{n=1}^{\infty },$
\ such that%
\begin{equation*}
\rho ^{n},\;(\rho ^{n})^{2}\rightharpoonup \rho ,\;\beta ^{\ast }\qquad 
\text{$\star $-weakly in $L^{\infty }(\Omega _{T})$}~\text{and}~\text{weakly
in $L^{2}(\Omega _{T}).$}
\end{equation*}%
Also, we have%
\begin{equation*}
\rho ^{o,n},\;(\rho ^{o,n})^{2}\rightharpoonup \rho ^{o},\;\beta ^{o,\ast
}\qquad \text{$\star $-weakly in }L^{\infty }(\Gamma _{T};\mu ^{+})\text{$.$}
\end{equation*}%
By construction the pair $(\rho ^{n},\mathbf{v}^{n})$ fulfills equality %
\eqref{1}, that yields the estimate 
\begin{equation}
\Bigl|\iint\limits_{\Omega _{T}}\rho ^{n}\,\phi _{t}\,d\mathbf{x}dt\Bigl|%
\leqslant C\,\Vert \phi \Vert _{L^{2}(0,T;H_{0}^{1}(\Omega ))},  \label{r}
\end{equation}%
which is valid for an arbitrary function $\phi \in C_{0}^{\infty
}(0,T;H_{0}^{1}(\Omega ))$. Also, applying Lemma \ref{prop22} the functions $%
\alpha =(\rho ^{n})^{2}$ solves system \eqref{2} (in the integral sense %
\eqref{3}), hence 
\begin{equation}
\Bigl|\iint\limits_{\Omega _{T}}(\rho ^{n})^{2}\,\phi _{t}\,d\mathbf{x}dt%
\Bigl|\leqslant C\,\Vert \phi \Vert _{L^{2}(0,T;H_{0}^{1}(\Omega ))},\quad
\forall \phi \in C_{0}^{\infty }(0,T;H_{0}^{1}(\Omega )).  \label{rr}
\end{equation}%
From the compact embedding of $L^{2}(\Omega )$ in $\ H^{-1}(\Omega ),$ %
\eqref{r}, \eqref{rr}, and the well known compactness results of \cite{aubin}%
, \cite{simon}, we have 
\begin{equation}
\rho ^{n},\;(\rho ^{n})^{2}\rightarrow \rho ,\;\beta ^{\ast }\quad 
\mbox{
strongly  in }L^{2}\bigl(0,T;H^{-1}(\Omega )\bigr).  \label{c2}
\end{equation}%
Since $\mathbf{v}^{n}-\mathbf{v}_{b}\in L^{2}(0,T;\mathbf{H}_{0}^{1}(\Omega
))$ (see \eqref{vb}),\ due to \eqref{up}, (\ref{c2}), we obtain 
\begin{equation*}
\iint_{\Omega _{T}}\rho ^{n}\left( (\mathbf{v}^{n}-\mathbf{v}_{b})\cdot {%
\nabla }\right) \phi \;d\mathbf{x}dt\rightarrow \iint_{\Omega _{T}}\rho
\left( (\mathbf{v}-\mathbf{v}_{b})\cdot {\nabla }\right) \phi \;d\mathbf{x}dt
\end{equation*}%
for each test functions $\phi \in C_{c}^{\infty }((-\infty ,T)\times \mathbb{%
R}^{d})$, that is to say 
\begin{equation}
\iint_{\Omega _{T}}\rho ^{n}(\mathbf{v}^{n}\cdot {\nabla )}\phi \;d\mathbf{x}%
dt\rightarrow \iint_{\Omega _{T}}\rho (\mathbf{v}\cdot {\nabla )}\phi \;d%
\mathbf{x}dt.  \label{c01}
\end{equation}%
Similarly, applying the same idea, we deduce 
\begin{equation}
\iint_{\Omega _{T}}(\rho ^{n})^{2}(\mathbf{v}^{n}\cdot {\nabla )}\phi \;d%
\mathbf{x}dt\rightarrow \iint_{\Omega _{T}}\beta ^{\ast }(\mathbf{v}\cdot {%
\nabla )}\phi \;d\mathbf{x}dt.  \label{c02}
\end{equation}

Let us pass to the limit $n\rightarrow \infty $ in equalities \eqref{1} and %
\eqref{3}, written for $(\rho ^{n},\rho ^{o,n}),$ $\mathbf{v}^{n}$ and $%
\alpha =((\rho ^{n})^{2},(\rho ^{o,n})^{2}),$ $\mathbf{v}^{n}$, \
respectively. Thus with the help of \eqref{c01}, \eqref{c02}, we obtain that
the triples $(\rho ,\rho ^{o}),$ $\mathbf{v}$ and $\alpha =(\beta ^{\ast
},\beta ^{o,\ast }),$ $\mathbf{v}$ fulfill equalities \eqref{1} and \eqref{3}%
. \ Applying Lemma \ref{prop22} to $(\rho ,\rho ^{o}),$ $\mathbf{v}$, we
have that, the triple $(\rho ^{2},(\rho ^{o})^{2}),$ $\mathbf{v}$ satisfies %
\eqref{3} too. In view of the uniqueness of solution for system \eqref{2},
we obtain 
\begin{equation*}
\beta ^{\ast }\equiv \rho ^{2}\quad \text{a.e. in }\Omega _{T}.
\end{equation*}%
Hence $\rho ^{n}\rightarrow \rho $ \ strongly in $L^{2}(\Omega _{T}).$
Therefore we derive that there exists a suitable sequence of $\left\{ (\rho
^{n},\nu ^{n})=P(\overline{\rho }^{n},\overline{\nu }^{n}):~(\overline{\rho }%
^{n},\overline{\nu }^{n})\in \mathcal{Z}\right\} _{n=0}^{\infty },$ such
that 
\begin{equation*}
(\rho ^{n},\nu ^{n})\rightarrow (\rho ,\nu )\quad \text{strongly in }%
\;L^{2}(\Omega _{T})^{2}\quad \text{and}\quad (\rho ,\nu )\in \mathcal{Z}.
\end{equation*}%
Consequently, $P$ is a compact operator on the set $\mathcal{Z}$.$\hfill
\;\blacksquare $

\bigskip

Now we show that

\begin{lemma}
\label{PPA21} The operator $P$ is continuous in the norm of $L^{2}(\Omega
_{T})^{2}.$
\end{lemma}

\textbf{Proof.} \ Let $\left\{ (\overline{\rho }^{n},\overline{\nu }^{n})\in 
\mathcal{Z}\right\} _{n=0}^{\infty }$ be a sequence converging to $(%
\overline{\rho },\overline{\nu })\in \mathcal{Z}$ in $L^{2}(\Omega
_{T})^{2}, $ that is 
\begin{equation}
||(\overline{\rho }^{n},\overline{\nu }^{n})-(\overline{\rho },\overline{\nu 
})||_{\text{$L^{2}(\Omega _{T})^{2}$}}\rightarrow 0,\quad n\rightarrow
\infty .  \label{op}
\end{equation}%
Let $(\rho ^{n},\nu ^{n})=P(\overline{\rho }^{n},\overline{\nu }^{n})\in 
\mathcal{Z},$ $\mathbf{v}^{n}\in L^{2}(0,T;\mathbf{V}^{1}(\Omega ))$\ \ and $%
(\rho ,\nu )=P(\overline{\rho },\overline{\nu })\in \mathcal{Z},$ $\ \mathbf{%
v}\in L^{2}(0,T;\mathbf{V}^{1}(\Omega ))$ \ be the solutions of \eqref{a}
and \eqref{b}, respectively. The triples $(\rho ^{n},\rho ^{o,n}),$ $(\nu
^{n},\nu ^{o,n}),$ $\mathbf{v}^{n}$ and $(\rho ,\rho ^{o}),$ $(\nu ,\nu
^{o}),$ $\mathbf{v}$ satisfy \eqref{EUU}. Let us consider $\mathbf{z}^{n}:=%
\mathbf{v}^{n}-\mathbf{v}$, $P^{n}:=p^{n}-p,$ which satisfy the system 
\begin{equation*}
\begin{cases}
-\mathrm{div}(\overline{\rho }^{n}\overline{\nu }^{n}\mathbf{Dz}^{n})+h(t,x,%
\overline{\rho }^{n}\overline{\nu }^{n})\,\mathbf{z}^{n}=-\nabla P^{n}+%
\mathbf{f}^{n},\quad div\mathbf{z}^{n}=0\quad \text{in }\Omega _{T}, \\ 
\mathbf{z}^{n}=\mathbf{0}\quad \text{on }\Gamma _{T},%
\end{cases}%
\end{equation*}%
with $\mathbf{f}^{n}:=(\overline{\rho }^{n}-\overline{\rho })\mathbf{g}+%
\mathrm{div}((\overline{\rho }^{n}\overline{\nu }^{n}-\overline{\rho }\ 
\overline{\nu })\mathbf{D}(\mathbf{v}))-s_{n}\mathbf{v}.$ \ \ Here, we
define 
\begin{equation*}
s_{n}:=h(t,\mathbf{x},\overline{\rho }^{n}\overline{\nu }^{n})-h(t,\mathbf{x}%
,\overline{\rho }\ \overline{\nu })\quad \text{and}\quad l_{n}:=|\overline{%
\rho }^{n}-\overline{\rho }\ |+|\overline{\nu }^{n}-\overline{\nu }|.
\end{equation*}%
Due to \eqref{reg3} and \eqref{EUU} the function$\ \mathbf{z}^{n}(t)$
satisfies for a.a. $t\in (0,T)$ the following 
\begin{align*}
\overline{\rho }^{0}\overline{\nu }^{0}\Vert \mathbf{Dz}^{n}\Vert _{\mathbf{L%
}^{2}(\Omega )}^{2}\leqslant \int_{\Omega }l_{n}|\mathbf{g}||\mathbf{z}%
^{n}|\,d\mathbf{x}+C\int_{\Omega }l_{n}|\mathbf{Dv}||\mathbf{Dz}^{n}|\,d%
\mathbf{x}& +\int_{\Omega }|s_{n}||\mathbf{v}||\mathbf{z}^{n}|\,d\mathbf{x}
\\
\leqslant C\left( ||l_{n}\mathbf{g}||_{\mathbf{L}^{q}(\Omega )}^{2}+||l_{n}%
\mathbf{Dv}||_{\mathbf{L}^{2}(\Omega )}^{2}\right) +C||s_{n}||_{\mathbf{L}%
^{s}(\Omega )}^{2}||\mathbf{v}||_{\mathbf{L}^{r}(\Omega )}^{2}& +\frac{%
\overline{\rho }^{0}\overline{\nu }^{0}}{4}\Vert \mathbf{Dz}^{n}\Vert _{%
\mathbf{L}^{2}(\Omega )}^{2},
\end{align*}%
where $\frac{1}{s}+\frac{2}{r}=1$ for $s$ given by \eqref{reg3}, and the
constant $C$ depend only on $\Omega ,$ $\overline{\nu }^{i},$ $\overline{%
\rho }^{i},$\ $i=0,1.$\ \ Here we have used Holder's and Young's
inequalities, inequality \eqref{ldd}, and the embedding \eqref{s}.
Therefore, we have%
\begin{equation*}
\Vert \mathbf{Dz}^{n}\Vert _{L^{2}(0,T;\mathbf{L}^{2}(\Omega
))}^{2}\leqslant C(I_{n}+J_{n}),
\end{equation*}%
where%
\begin{equation*}
I_{n}=||l_{n}\mathbf{g}||_{L^{2}(0,T;\mathbf{L}^{q}(\Omega ))}^{2}+||l_{n}%
\mathbf{Dv}||_{\mathbf{L}^{2}(\Omega _{T})}^{2},\quad
J_{n}=\int_{0}^{T}||s_{n}||_{\mathbf{L}^{s}(\Omega )}^{2}||\mathbf{v}||_{%
\mathbf{H}^{1}(\Omega )}^{2}dt.
\end{equation*}

\medskip Now, let us define the cut--off function $\phi _{L}(t,\mathbf{x}%
)=\min ~\{\phi (t,\mathbf{x}),L\}\;$ for $L>0.$ We have 
\begin{eqnarray*}
I_{n} &\leqslant &CL^{2}||(\overline{\rho }^{n},\overline{\nu }^{n})-(%
\overline{\rho },\overline{\nu })||_{\text{$L^{2}(\Omega _{T})^{2}$}%
}^{2}+C\int_{0}^{T}\left( \int_{\Omega }||\mathbf{g}|^{q}-|\mathbf{g}%
|_{L}^{q}|\,d\mathbf{x}\right) ^{2/q}dt \\
&&+C\iint_{\Omega _{T}}||\mathbf{Dv}|^{2}-|\mathbf{Dv}|_{L}^{2}|\,d\mathbf{x}%
dt.
\end{eqnarray*}%
Choosing $L:=||(\overline{\rho }^{n},\overline{\nu }^{n})-(\overline{\rho },%
\overline{\nu })||_{\text{$L^{2}(\Omega _{T})^{2}$}}^{-1/2},$ we derive \ $%
I_{n}\rightarrow 0$ as $n\rightarrow \infty .$

Also we have%
\begin{equation*}
J_{n}\rightarrow 0\quad \text{as \ }n\rightarrow \infty ,
\end{equation*}%
which can be proved by contradiction. Indeed, let us assume that there
exists a subsequence of $\left\{ J_{n}\right\} _{n=0}^{\infty }$ $,$\ which
is not convergent to zero as \ $n\rightarrow \infty :$ 
\begin{equation}
J_{n^{\prime }}\geqslant const>0,\ \ \forall n^{\prime }.  \label{01}
\end{equation}%
But, using \eqref{op}, there exists a subsequence $\left\{ s_{n^{^{\prime
\prime }}}\right\} _{n^{^{\prime \prime }}=0}^{\infty }$ of $\left\{
s_{n^{\prime }}\right\} _{n^{\prime }=0}^{\infty },$ such that $s_{n^{\prime
\prime }}\rightarrow 0$ a.e. in $\Omega _{T},$ since $h$ is a Carath\'{e}%
odory function (see \eqref{reg3}).\ Applying the uniform $L^{\infty }$%
--boundness of $\{(\overline{\rho }^{n^{\prime \prime }},$ $\overline{\nu }%
^{n^{\prime \prime }})\}_{n^{\prime \prime }=0}^{\infty },$ $(\overline{\rho 
},\overline{\nu })$ and \eqref{reg3}, we get 
\begin{equation*}
\int_{0}^{T}||s_{n^{\prime \prime }}||_{\mathbf{L}^{s}(\Omega
)}(t)dt\rightarrow 0
\end{equation*}%
by the dominated convergence theorem. Moreover Fubini's theorem, the
relation between the convergence in $L^{1}$--norm and the pointwise
convergence give the existence of a suitable subsequence $\left\{
s_{n^{^{\prime \prime \prime }}}\right\} _{n^{^{\prime \prime \prime
}}=0}^{\infty }$ of $\left\{ s_{n^{\prime \prime }}\right\} _{n^{\prime
\prime }=0}^{\infty },$ such that $||s_{n^{\prime \prime \prime }}||_{%
\mathbf{L}^{s}(\Omega )}(t)\rightarrow 0$ a.e. in $(0,T).$ \ Therefore the
uniform $L^{\infty }$--boundness of $\left\{ ||s_{n^{\prime \prime \prime
}}||_{\mathbf{L}^{s}(\Omega )}(t)\right\} _{n^{\prime \prime }=0}^{\infty }$
\ (see \eqref{reg3})\ and the dominated convergence theorem imply%
\begin{equation*}
J_{n^{\prime \prime \prime }}\rightarrow 0\text{\quad as }n^{\prime \prime
\prime }\rightarrow \infty ,
\end{equation*}%
which is a contradiction with \eqref{01}. Therefore we conclude%
\begin{equation}
\Vert \mathbf{v}^{n}-\mathbf{v}\Vert _{L^{2}(0,T;\mathbf{V}^{1}(\Omega
))}\rightarrow 0\quad \text{as }n\rightarrow 0.  \label{o2}
\end{equation}

Now, using \eqref{o2}, we can prove%
\begin{equation*}
r_{n}:=||(\rho ^{n},\nu ^{n})-(\rho ,\nu )||_{\text{$L^{2}(\Omega _{T})^{2}$}%
}\mathop{\rightarrow}\limits_{n\rightarrow \infty }0
\end{equation*}%
by contradiction.\ \ Let us assume that there exists a subsequence of $%
\left\{ r_{n}\right\} _{n=0}^{\infty },$ \ which is not convergent to zero
as \ $n\rightarrow \infty :$ 
\begin{equation}
r_{n^{\prime }}\geqslant const>0,\ \ \forall n^{\prime }.  \label{02}
\end{equation}%
Although, there exists a subsequence%
\begin{align*}
& \rho ^{n^{\prime \prime }},\;\nu ^{n^{\prime \prime }}\rightharpoonup 
\widetilde{\rho },\;\widetilde{\nu }\qquad \star \text{-weakly in }%
\;L^{\infty }(\Omega _{T})\quad \text{and}\quad \text{weakly in }%
L^{2}(\Omega _{T}), \\
& \rho ^{o,n^{\prime \prime }},\;\nu ^{o,n^{\prime \prime }}\rightharpoonup 
\widetilde{\rho }^{o},\;\widetilde{\nu }^{o}\qquad \text{$\star $-weakly in }%
L^{\infty }(\Gamma _{T};\mu ^{+}), \\
& \mathbf{v}^{n^{\prime \prime }}\rightarrow \mathbf{v}\qquad \text{strongly
in }\;L^{2}(0,T;\mathbf{V}^{1}(\Omega )),
\end{align*}%
where $(\widetilde{\rho },\widetilde{\rho }^{o}),$ $(\widetilde{\nu },%
\widetilde{\nu }^{o}),$ $\mathbf{v}$ satisfy \eqref{EUU} and solve systems %
\eqref{b}, respectively. In view of the uniqueness of solution for the two
transport systems in \eqref{b} (uniqueness is an important point here), for
the given $\mathbf{v,}$ we conclude that $\widetilde{\rho }\equiv \rho ,\;%
\widetilde{\nu }\equiv \nu .$ Now we can argue as in Lemma \ref{PPA2} and
show the existence of a subsequence%
\begin{equation*}
\rho ^{n^{\prime \prime \prime }},\;\nu ^{n^{\prime \prime \prime
}}\rightarrow \rho ,\;\nu \quad \text{strongly in}\;L^{2}(\Omega _{T}),
\end{equation*}%
which is a contradiction with our assumption \eqref{02}. Consequently, we
have shown the continuity of $P$. $\blacksquare $\newline

\bigskip The thesis of Theorem \ref{theorem} follows combining Lemmas \ref%
{PPA}, Lemma \ref{PPA2} and Lemma \ref{PPA21}.


\subsection{Two--Phase Filtration}

\label{two} 

In this section we show that the system $\mathbf{GMP}$ describes the motion
of \textit{immiscible} fluids.

\medskip We assume that the fluids\ $"1"$ and $"2"$ occupy two measurable
disjoint sets $\Omega _{1}$ and $\Omega _{2}$ at the initial moment $t=0,$
such that $\Omega _{1}\cup \Omega _{2}=\Omega .$ Moreover the fluids\ $"1"$
and $"2"$ enter inside of the domain $\Omega $ \ through two separated
boundary zones $\Gamma _{1}^{-}$ \ and $\Gamma _{2}^{-},$ such that $\Gamma
_{1}^{-}\cup \Gamma _{2}^{-}=\Gamma _{T}^{-}.$\ \ Let the data $\rho
_{0},\rho _{b}$ and \ $\nu _{0},\nu _{b}$ satisfy the natural assumptions
for $i=1,2:$ 
\begin{equation}
(\rho _{0},\nu _{0}):=(\rho _{i,0},\nu _{i,0})\quad \text{in }\Omega
_{i}\quad \text{and}\quad (\rho _{b},\nu _{b}):=(\rho _{i,b},\nu
_{i,b})\quad \text{on }\Gamma _{i}^{-}  \label{reg4}
\end{equation}%
with 
\begin{eqnarray}
\rho _{i,0} &\in &[\overline{\rho }_{i}^{0},\overline{\rho }_{i}^{1}],\quad
\nu _{i,0}\in \lbrack \overline{\nu }_{i}^{0},\overline{\nu }_{i}^{1}]\quad 
\text{a.e. in }\Omega _{i}(0),  \notag \\
\rho _{i,b} &\in &[\overline{\rho }_{i}^{0},\overline{\rho }_{i}^{1}],\quad
\nu _{i,b}\in \lbrack \overline{\nu }_{i}^{0},\overline{\nu }_{i}^{1}]\quad 
\text{a.e. on }\Gamma _{i}^{-},  \notag
\end{eqnarray}%
for positive real numbers $\overline{\rho }_{1}^{0}\leqslant \overline{\rho }%
_{1}^{1}<$ $\overline{\rho }_{2}^{0}\leqslant \overline{\rho }_{2}^{1}$ $\ $%
and $\overline{\nu }_{1}^{0}\leqslant \overline{\nu }_{1}^{1}<$ $\overline{%
\nu }_{2}^{0}\leqslant \overline{\nu }_{2}^{1}.$

\bigskip In the following, we demonstrate that the process of motion is 
\textit{immiscible.}

\begin{corollary}
\label{p2} \label{prop233} Let the data $\rho _{0},\nu _{0}$, \ $\rho
_{b},\nu _{b}$, and $\mathbf{b}$ fulfill the assumptions of Theorem \ref%
{theorem}. If $\rho _{0},\nu _{0}$, \ $\rho _{b},\nu _{b}$ satisfy the
additional assumptions \eqref{reg4}, then the weak solution $(\rho ,\rho
^{o}),$ $(\nu ,\nu ^{o}),$ $\mathbf{v}$ of $\mathbf{GMP}$ describes the
behaviour of two \underline{\textit{immiscible}} \ fluids $"1"$ and $"2"$ in
the porous media, that is, there exist two measurable disjoint sets $F_{1},$%
\ $F_{2}$ with $F_{1}\cup F_{2}=\Omega _{T}$ and two measurable disjoint
boundary zones $\Gamma _{1}^{+},$ $\Gamma _{2}^{+}$ with $\Gamma
_{1}^{+}\cup \Gamma _{2}^{+}=\Gamma _{T}^{+},$ such that%
\begin{eqnarray*}
\rho &\in &[\overline{\rho }_{i}^{0},\overline{\rho }_{i}^{1}],\quad \nu \in
\lbrack \overline{\nu }_{i}^{0},\overline{\nu }_{i}^{1}]\quad \text{a.e. in }%
F_{i}, \\
\rho ^{o} &\in &[\overline{\rho }_{i}^{0},\overline{\rho }_{i}^{1}],\quad
\nu ^{o}\in \lbrack \overline{\nu }_{i}^{0},\overline{\nu }_{i}^{1}]\quad 
\text{a.e. on }\Gamma _{i}^{+}
\end{eqnarray*}%
for $i=1,2$, respectively.
\end{corollary}

\textbf{Proof. \ } 1. First, let us take an arbitrary non--negative  function $%
\beta \in C^{1}(\mathbb{R})$ with ${\mathrm{supp}(\beta )}= \overline{%
\mathrm{\mathbb{R}} \setminus \cup _{i=1,2}[\overline{\rho }_{i}^{0}, 
\overline{\rho }_{i}^{1}]}.$ By the same arguments as in proof of Lemma \ref%
{prop23}, we obtain that $\alpha =(\beta (\rho ),\beta (\rho ^{o}))$ is the
unique weak solution of system \eqref{2} and $\alpha \equiv 0$. Therefore,
the functions $\rho ,$ $\rho ^{o}$\ take values just in $\cup _{i=1,2}[%
\overline{\rho }_{i}^{0},\overline{\rho }_{i}^{1}]$.

\medskip 2. Now, taking an arbitrary non--negative function $\beta _{1}\in
C^{1}(\mathbb{R})$ with $\mathrm{supp}(\beta _{1})= [\overline{\rho }%
_{1}^{0},\overline{\rho }_{1}^{1}],$ we consider the measurable sets 
\begin{eqnarray*}
F_{1}^{\rho } &:&=\left\{ (t,\mathbf{x})\in \Omega _{T}:~\beta _{1}(\rho (t,%
\mathbf{x}))>0\right\} ,\quad F_{2}^{\rho }:=\Omega _{T}\backslash
F_{1}^{\rho }, \\
\Gamma _{1}^{\rho ^{o}} &:&=\left\{ (t,\mathbf{x})\in \Gamma _{T}^{+}:~\beta
_{1}(\rho ^{o}(t,\mathbf{x}))>0\right\} ,\quad \Gamma _{2}^{\rho
^{o}}:=\Gamma _{T}^{+}\backslash \Gamma _{1}^{\rho ^{0}},
\end{eqnarray*}%
which satisfy 
\begin{equation*}
\rho \in \lbrack \overline{\rho }_{i}^{0},\overline{\rho }_{i}^{1}]\quad 
\text{a.e. in }F_{i}^{\rho }\quad \text{and}\quad \rho ^{o}\in \lbrack 
\overline{\rho }_{i}^{0},\overline{\rho }_{i}^{1}]\quad \text{a.e. on }%
\Gamma _{i}^{\rho ^{o}},\quad i=1,2.
\end{equation*}%
Analogously, we can also define the sets $F_{1}^{\nu }:=\left\{ (t,\mathbf{x}%
)\in \Omega _{T}:~\beta _{2}(\nu (t,\mathbf{x}))>0\right\} ,$ $F_{2}^{\nu
}:=\Omega _{T}\backslash F_{1}^{\nu }$ and $\Gamma _{1}^{\nu ^{o}}:=\left\{
(t,\mathbf{x})\in \Gamma _{T}^{+}:~\beta _{2}(\nu ^{o}(t,\mathbf{x}%
))>0\right\} ,$ $\Gamma _{2}^{\nu ^{o}}:=\Gamma _{T}^{+}\backslash \Gamma
_{1}^{\nu ^{o}},$ such that 
\begin{equation*}
\nu \in \lbrack \overline{\nu }_{i}^{0},\overline{\nu }_{i}^{1}]\quad \text{%
a.e. in }F_{i}^{\nu }\quad \text{and}\quad \nu ^{o}\in \lbrack \overline{\nu 
}_{i}^{0},\overline{\nu }_{i}^{1}]\quad \text{a.e. on }\Gamma _{i}^{\nu
^{o}},\quad i=1,2,
\end{equation*}%
for an arbitrary non--negative function $\beta _{2}\in C^{1}(\mathbb{R})$
with $\mathrm{supp}(\beta _{2})= [\overline{\nu }_{1}^{0},\overline{\nu }%
_{1}^{1}]. $

\bigskip By Proposition 4.2 in \cite{boyer}, we see that the $\alpha =(\beta
_{1}(\rho )(1-\beta _{2}(\nu )),~\beta _{1}(\rho ^{o})(1-\beta _{2}(\nu
^{o})))$ and $\alpha =((1-\beta _{1}(\rho ))\beta _{2}(\nu ),~(1-\beta
_{1}(\rho ^{o}))\beta _{2}(\nu ^{o}))$\ satisfy the following system (in the
distributional sense)%
\begin{eqnarray*}
\partial _{t}\alpha + \mathrm{div} (\mathbf{v}\alpha ) &=&0\quad \text{in }%
\Omega _{T}, \\
\alpha |_{t=0} &=&0\quad \text{in }\Omega \quad \text{and}\quad \alpha
|_{\Gamma _{T}^{-}}=0.
\end{eqnarray*}%
Moreover, from the uniqueness of solution of this system, we have $\alpha
\equiv 0$. Therefore, we have 
\begin{equation*}
F_{i}^{\rho }\equiv F_{i}^{\nu }=:F_{i}\quad \text{and}\quad \Gamma
_{i}^{\rho ^{o}}\equiv \Gamma _{i}^{\nu ^{o}}=:\Gamma _{i}^{+},\quad i=1,2.
\end{equation*}%
The sets $F_{i}$ and $\Gamma _{i}^{+},\ i=1,2,$ defined above, do not depend
on the choose of $\beta _{1},\beta _{2}$. Indeed, it is enough to apply one
more time Proposition 4.2 in \cite{boyer}.$\hfill \blacksquare $

\bigskip

Finally, as a particular result we have the following

\begin{corollary}
\label{prop2333} Let the data $\rho _{0},\nu _{0}$, \ $\rho _{b},\nu _{b}$,
and $\mathbf{b}$ fulfill the assumptions of Lemma \ref{prop233}. Let $\rho
_{i}:=\overline{\rho }_{i}^{0}\equiv \overline{\rho }_{i}^{1}$ and $\nu
_{i}:=\overline{\nu }_{i}^{0}\equiv \overline{\nu }_{i}^{1},$ $\ i=1,2,$\
then there exist measurable disjoint sets $F_{1},$\ $F_{2}$ with $F_{1}\cup
F_{2}=\Omega _{T}$ \ and two measurable disjoint boundary zones $\Gamma
_{1}^{+},$ $\Gamma _{2}^{+}$ with $\Gamma _{1}^{+}\cup \Gamma
_{2}^{+}=\Gamma _{T}^{+},$ such that 
\begin{eqnarray*}
(\rho ,\nu ) &=&(\rho _{i},\nu _{i})\quad \text{a.e. in }F_{i}, \\
(\rho ^{o},\nu ^{o}) &=&(\rho _{i},\nu _{i})\quad \text{a.e. on }\Gamma
_{i}^{+}\quad \text{for }i=1,2.
\end{eqnarray*}
\end{corollary}

\bigskip \textbf{Open problems:}

1) Let us remark that the uniqueness result for the problem $\mathbf{GMP}$%
\textbf{\ }is not shown.

\medskip 2) The investigation of the regularity of the interface $S(t)$,
which separates the two fluids during the motion, is an interesting open
problem. Let us point that such problem was studied in the article \cite%
{ACDCFG} for the interface $S(t)$ between two incompressible 2--D fluids,
where the evolution equation for $S(t)$ was obtained from \underline{Darcy's
law}. \ 


\section{Acknowledgements}


We thank the anonymous referee for providing very constructive comments and
help us improve the contents of this paper.

\smallskip The work of N.V. Chemetov was supported through the project
POCTI/ISFL / 209 of CMAF/UL and the project PTDC/MAT/110613/2009 of
FFC/FC/UL. 
Wladimir Neves is partially supported by
Pronex-FAPERJ through the grant E-26/ 110.560/2010 entitled \textsl{"Nonlinear Partial
Differential Equations"}, also by CNPq through the grants 
484529/2013-7, 308652/2013-4.



\begin{thebibliography}{99}

\bibitem{m}  Alarc\'on E.A.,  Del Sol M.M., I\'orio Junior R.J., 
\emph{On the Cauchy problem associated to the Brinkman flow in $\mathbb{R}^n$.}
Applicable Analysis and Discrete Mathematics, \textbf{6}, n. 2 (2012), 214--237.

\bibitem{al} Allaire G., \emph{Homogenization of the Navier--Stokes equations
in open sets perforated with tiny holes. I: Abstract framework, a volume
distribution of holes.} Arch. Ration.Mech. Anal., \textbf{113}, n. 3 (1991),
209--259.

\bibitem{ambrose} Ambrose D., \emph{Well--posedness of two--phase Hele--Shaw
Flow without surface tension.} Euro. Jnl. of Applied Mathematics, \textbf{15}
(2004), 597--607.



\bibitem{A1} Antontsev S.N., Meirmanov A., Yurinsky B.V., \emph{ A free-boundary problem for Stokes equations:
classical solutions.} Interfaces and Free Boundaries, \textbf{2} (2000), 413-424.


\bibitem{RJAREC} Atkin R.J., Craine R.E., \emph{Continuum theories of
mixtures: Applications}. Journal of the Institute of Mathematics and its
Applications, \textbf{17} (1976), 153--207.

\bibitem{aubin} Aubin J.--P., \emph{Un theoreme de compacite. }\textit{C. R.
Acad. Sci. = Paris,} \textbf{256} (1963), 5042--5044.

\bibitem{aur} Auriault J.--L., \emph{On the Domain of Validity of Brinkman's
Equation. }Transport Porous Media, \textbf{79} (2009), 215--223




\bibitem{gr} Bosia S., Conti M.,  Grasselli M., \emph{On the Cahn--Hillard--Brinkman system,} submitted.


\bibitem{boyer} Boyer F., \emph{Trace theorems and spatial continuity
properties for the solutions of the transport equation.} Differential and
integral equations, \textbf{18}, n. 8 (2005), 891--934.

\bibitem{brin} Brinkman H.C., \emph{A calculation of the viscouse force
exerted by a flowing fluid on a dense swarm of particles. } Appl. Sci. Res., 
\textbf{A1} (1947), 27--34.

\bibitem{cat} Cattabriga L., \emph{Su un problema al contorno relativo al
sistema di equazioni di Stokes (Italian).} Rend. Sem. Mat. Univ. Padova, 
\textbf{31} (1961), 308--340.


\bibitem{ch} Chand R., Rana G.C., \emph{On the onset of thermal convection in rotating nanofluid 
layer saturating a Darcy--Brinkman porous medium.} Int. J. of Heat and Mass Transfer, \textbf{55} (2012), 5417-- 5424.

\bibitem{NCWN1} Chemetov N., Neves W., \emph{The generalized
Buckley--Leverett system. Solvability}. Arch. Ration. Mech. Analysis, \textbf{%
208} (2013), 1--24.

\bibitem{const} Constantin P., Foias C., \emph{Navier--Stokes equations.}
Chicago Lectures in Mathematics, The University of Chicago Press, Ltd.,
London, 1988.

\bibitem{PCDCFGRMS} Constantin P., C\'{o}rdoba D., Gancedo F., Strain R.M., 
\emph{On the global existence for the Muskat problem}. J. Eur. Math. Soc., 
\textbf{15} (2013), 201--227.

\bibitem{cor} Cordoba D., Faraco D., Gancedo F., \emph{Lack of Uniqueness
for Weak Solutions of the Incompressible Porous Media Equation,} Arch.
Rational Mech. Anal., \textbf{200} (2011), 725--746.

\bibitem{ACDCFG} C\'{o}rdoba A., C\'{o}rdoba D., Gancedo F., \emph{Interface
evolution: the Hele--Shaw and Muskat problems}. Ann. of Math., \textbf{173}
(2011), 477--542.


\bibitem{den} Denisova I.V., Solonnikov V.A., 
\emph{Global solvability of the problem of the motion of two incompressible capillary fluids in a container.} (Russian) 
Zap. Nauchn. Sem. S.--Peterburg. Otdel. Mat. Inst. Steklov. (POMI) \textbf{397} (2011), 
Kraevye Zadachi Mat Fiziki i Smezhnye Voprosy Teorii Funktsii. 42, 20--52, 172;
 English translation in J. Math. Sci. (N. Y.) \textbf{185}, n. 5 (2012), 668-–686.

\bibitem{JEBVM} Escher J., Matioc B.--V., \emph{On the parabolicity of the
Muskat Problem: Well--posedness, fingering, and stability results.} 
 Z. Anal. Anwend., \textbf{30}, n. 2 (2011), 193--218.

\bibitem{EG} Evans L.C., \emph{Partial differential equations.} Graduate
Studies in Mathematics, Volume 19, AMS, Providence, Rhode Island, 1998.


\bibitem{gir} Girault V., Raviart P.--A.,\emph{ Finite element methods for the
Navier--Stokes equations. Theory and Algorithms.} Springer--Verlag, Berlin
Heidelberg New York Tokyo, 1986.

\bibitem{gir2} Girault V., Kanschat G., Riviere B., 
\emph{On the Coupling of Incompressible Stokes or Navier–Stokes and Darcy Flows Through Porous Media.} 
Modelling and Simulation in Fluid Dynamics in Porous Media. Springer Proceedings in Mathematics and Statistics, 
\textbf{28} (2013), 1--25.

\bibitem{i} Ingram, R., 
\emph{ A Mixed Finite Element Approximation of Stokes-Brinkman and NS--Brinkman Equation for Non--Darcian Flows.}
 SIAM J. Numer. Anal.,  \textbf{49}, 2 (2011), 491-520.


\bibitem{k} Kelliher J.P., Temam R., Wang X., \emph{ Boundary layer associated with 
the Darcy-–Brinkman-–Boussinesq model for convection in porous media.}
Physica D,  \textbf{240}, 7 (2011), 619-–628.



\bibitem{kl} Khaled A.R.A., Vafai K., \emph{The role of porous media in modeling flow and heat transfer in
biological tissues.} Int. J. Heat Mass Transfer,  \textbf{46} (2003),  4989--5003.




\bibitem{KLLS} Krotkiewski M., Ligaarden I.S., Lie K.--A., Schmid D.W., 
\emph{On the Importance of the Stokes--Brinkman Equations for Computing
Effective Permeability in Carbonate Karst Reservoirs}.
Commun. Comput. Phys., (2011), 1--18.




\bibitem{L1} Layton W.J.,  Schieweck F., Yotov I., \emph{Coupling fluid 
flow with porous media flow.} SIAM J. Numer. Anal., \textbf{40} (2003),  2195--2218.




\bibitem{MM} Massoudi M., \emph{Constitutive relations for the interaction
force in multicomponent particulate flows.} Int. Journal of Non--linear
Mechanics, \textbf{38} (2001), 313--336.

\bibitem{mar} Maru\v{s}ic--Paloka E., Pazanin I., Maru\v{s}ic S., \emph{%
Comparison between Darcy and Brinkman laws in a fracture.} Applied
Mathematics and Computation, \textbf{218} (2012), 7538--7545.

\bibitem{mo} Mosthafa K., Babera K., Flemischa B., 
Helmiga R., Leijnseb A., Rybakc I., Wohlmuthd B., 
\emph{A coupling concept for two--phase compositional
porous--medium and single--phase compositional free}.  
SRC SimTech, Univeritat Stutgart,  Preprint Series, Issue No. 2011--11, 1--23.


\bibitem{MUSKAT} Muskat M., \emph{Two fluid system in porous media. The
encroachment of water into oil sand.} Physics, \textbf{5} (1934), 250--264.

\bibitem{DANAB} Nield D.A., Bejan A., \emph{Convection in Porous Media.} Springer, 4th. Ed.,
2013.

\bibitem{KRRL} Rajagopal K.R., \emph{On a hierarchy of approximate models
for flows ~of incompressible fluids through porous solids.} Math. Models and
Methods in Applied Sciences, \textbf{17}, n. 2 (2012) 215--252.

\bibitem{KRRLT} Rajagopal K.R., Tao L., \emph{Mechanics of mixtures.} World
Scientific Publishing Co. Pte. Ltd., 1995.

\bibitem{MS} Sahimi M., \emph{Flow and Transport in Porous Media and
Fractured Rock.}  Wiley-VCH Verlag GmbH \& Co. KGaA, 2th. Ed., 1998.

\bibitem{Scheidegger} Scheidegger A.E., \emph{The Physics of Flow Through
Porous Media}, 3rd ed, University of Toronto Press, Toronto, 1974.

\bibitem{siegel} Siegel M., Caflisch R., Howison S., \emph{Global existence,
singular solutions, and ill--posedness for the Muskat problem.} \ Comm. Pure
and Appl. Math., \textbf{57} (2004), 1374--1411.

\bibitem{simon} Simon J., \emph{Compact sets in the space} $L\sp p(0,T;B)$. 
\textit{Ann. Mat. Pura Appl.} IV. Ser. \textbf{146} (1987), 65--96.


\bibitem{sze} Sz\'{e}kelyhidi Jr. \ L., \emph{Relaxation of the
incompressible porous media equation, }See pre--print: \
arxiv.org/pdf/1102.2597.


\bibitem{CTRT} Trusdell C., Toupin R., \emph{The classical field theories}.
Handbuch der Physik (Ed. S. Fl\"{u}gge) Vol. III/1 p.226, Springer--Verlag,
Berlin, 1960.

\bibitem{yi} Yi F., \emph{Global classical solution of Muskat free boundary
problem.} J. Math. Anal. Appl., \textbf{288} (2003), 442 -- 461.

\bibitem{y} Yosida K., \emph{Functional Analysis}. Springer--Verlag, Berlin,
1978.
\end{thebibliography}
\end{document}